\documentclass{amsart}
\usepackage{verbatim,amssymb,amsfonts,amscd}
\usepackage[all]{xy}
\numberwithin{equation}{section}

\let\cal\mathcal
\def\Ascr{{\cal A}}
\def\Bscr{{\cal B}}
\def\Cscr{{\cal C}}
\def\Dscr{{\cal D}}
\def\Escr{{\cal E}}
\def\Fscr{{\cal F}}

\def\Hscr{{\cal H}}

\def\Jscr{{\cal J}}

\def\Lscr{{\cal L}}
\def\Mscr{{\cal M}}

\def\Oscr{{\cal O}}
\def\Pscr{{\cal P}}
\def\Qscr{{\cal Q}}

\def\Sscr{{\cal S}}
\def\Tscr{{\cal T}}
\def\Uscr{{\cal U}}
\def\Vscr{{\cal V}}

\let\blb\mathbb

\def \PP{{\blb P}}

\def \ZZ{{\blb Z}}

\def \SS{{\blb S}}

\def\id{\text{id}}
\def\Id{\operatorname{id}}
\def\pr{\mathop{\text{pr}}\nolimits}

\let\at\ast

\def\Bimod{\operatorname{Bimod}}

\def\Mod{\operatorname{Mod}}
\def\mod{\operatorname{mod}}
\def\Gr{\operatorname{Gr}}

\def\QGr{\operatorname{QGr}}
\def\qgr{\operatorname{qgr}}

\def\gr{\operatorname{gr}}

\def\length{\mathop{\text{length}}}

\def\Supp{\mathop{\text{\upshape Supp}}}

\def\Qch{\mathop{\text{Qch}}}
\def\coh{\mathop{\text{\upshape{coh}}}}

\def\gr{\operatorname {gr}}
\def\Spec{\operatorname {Spec}}

\def\Ext{\operatorname {Ext}}
\def\Hom{\operatorname {Hom}}

\def\RHom{\operatorname {RHom}}
\def\uRHom{\operatorname {R\mathcal{H}\mathit{om}}}

\def\coker{\operatorname {coker}}
\def\ker{\operatorname {ker}}

\def\Tor{\operatorname {Tor}}

\def\id{{\operatorname {id}}}

\def\rk{\operatorname {rk}}

\def\Pic{\operatorname {Pic}}

\def\r{\rightarrow}

\DeclareMathOperator{\Proj}{Proj}

\DeclareMathOperator{\Tors}{Tors}
\DeclareMathOperator{\tors}{tors}

\DeclareMathOperator{\HHom}{\mathcal{H}\mathit{om}}

\DeclareMathOperator{\Pro}{Pro}

\let\dirlim\injlim
\let\invlim\projlim

\newtheorem{lemma}{Lemma}[section]

\newtheorem{theorem}[lemma]{Theorem}

\newtheorem{lemmas}{Lemma}[subsection]
\newtheorem{propositions}[lemmas]{Proposition}
\newtheorem{theorems}[lemmas]{Theorem}
\newtheorem{corollarys}[lemmas]{Corollary}

\theoremstyle{definition}

\newtheorem{definition}[lemma]{Definition}

\newtheorem{examples}[lemmas]{Example}
\newtheorem{definitions}[lemmas]{Definition}

{

\newtheorem{step}{Step}

}

\theoremstyle{remark}

\newtheorem{remarks}[lemmas]{Remark}

\newdimen\uboxsep \uboxsep=1ex
\def\uboxn#1{\vtop to 0pt{\hrule height 0pt depth 0pt\vskip\uboxsep
\hbox to 0pt{\hss #1\hss}\vss}}

\def\uboxs#1{\vbox to 0pt{\vss\hbox to 0pt{\hss #1\hss}
\vskip\uboxsep\hrule height 0pt depth 0pt}}

\date{\today}

\def\shbimod{\operatorname{shbimod}}
\def\ShBimod{\operatorname{ShBimod}}
\def\Points{\mathop{\mathrm{Points}}\nolimits}

\def\Sch{\mathop{\mathrm{Sch}}\nolimits}

\author{M. Van den Bergh}
 \address{Departement WNI,  Universiteit Hasselt, 
 3590 Diepenbeek, Belgium.}
  \email{michel.vandenbergh@uhasselt.be}
\thanks{The  author is a senior researcher at the FWO}
\date{\today}
\title{Non-commutative $\PP^1$-bundles over commutative schemes}
\makeatletter
\relax 
\@writefile{toc}{\contentsline {section}{\tocsection {}{1}{Introduction}}{1}}
\@writefile{toc}{\contentsline {section}{\tocsection {}{2}{Completion of abelian categories}}{3}}
\@writefile{toc}{\contentsline {subsection}{\tocsubsection {}{2.1}{Base extension}}{3}}
\@writefile{toc}{\contentsline {subsection}{\tocsubsection {}{2.2}{Completion of noetherian abelian categories}}{4}}
\@writefile{toc}{\contentsline {subsection}{\tocsubsection {}{2.3}{Functors}}{8}}
\@writefile{toc}{\contentsline {subsection}{\tocsubsection {}{2.4}{Formal flatness}}{10}}
\@writefile{toc}{\contentsline {subsection}{\tocsubsection {}{2.5}{$\operatorname  {Ext}$-groups}}{12}}
\@writefile{toc}{\contentsline {subsection}{\tocsubsection {}{2.6}{The complete derived category}}{15}}
\@writefile{toc}{\contentsline {section}{\tocsection {}{3}{Formal deformations of abelian categories}}{16}}
\@writefile{toc}{\contentsline {section}{\tocsection {}{4}{Ampleness}}{16}}
\@writefile{toc}{\contentsline {section}{\tocsection {}{5}{Lifting and base change}}{19}}
\@writefile{toc}{\contentsline {section}{\tocsection {}{}{References}}{19}}
\makeatother
\begin{document}
\begin{abstract}
In this paper we develop the theory of non-commutative $\PP^1$-bundles
over commutative (smooth) schemes. Such non-commutative
$\PP^1$-bundles occur in the theory of $D$-modules 
 but our definition is more general. We can show that every
non-commutative deformation of a Hirzebruch surface is given by a
non-commutative $\PP^1$-bundle over~$\PP^1$ in our sense.
\end{abstract}
\maketitle
\tableofcontents
\section{Introduction}
In this paper we develop the theory of non-commutative $\PP^1$-bundles
over commutative (smooth) schemes. Such non-commutative
$\PP^1$-bundles occur in the theory of $D$-modules (see
\cite{NevinsBenZvi}) but our definition is more general. The extra
generality is needed to cover basic examples in
non-commutative algebraic geometry \cite{VdB11}. As an
indication that our definition is the ``right one'' we present a proof that every
non-commutative deformation of a Hirzebruch surface is given by a
non-commutative $\PP^1$-bundle over~$\PP^1$ (see below).

\medskip

Let us explain our definition. Assume that $X$ is a scheme 
of finite type over a field $k$. Following \cite{VdB11} and later
\cite{Pat2,Pat1} we define a $\shbimod(X-X)$ as the category of coherent
$\Oscr_{X\times X}$ modules whose support is finite over $X$ on the
left and right. We call the elements of $\shbimod(X-X)$
``sheaf-bimodules'' to distinguish them from the somewhat more general
bimodules which were introduced in \cite{VdB19}. The category of
coherent sheaves on $X$ may be identified with the objects in
$\shbimod(X-X)$ supported on the diagonal. 

Convolution makes $\shbimod(X-X)$ into a monoidal category so we
may define a ``$\ZZ$-graded sheaf-algebra'' on $X$ to be a graded
algebra object  in $\shbimod(X-X)$. 
If
$\Ascr$ is a graded sheaf-algebra then we may define a category $\Gr(\Ascr)$
of graded $\Ascr$-modules.  Following \cite{AZ} we define $\QGr(\Ascr)$ as
$\Gr(\Ascr)$ divided by the modules which are direct limits of
right bounded ones. 

A first approximative approach to  non-commutative $\PP^1$-bundles on $X$,
advocated in \cite{Pat2,Pat1,VdB11}, is to consider abelian categories of the
form $\QGr(\Ascr)$ where $\Ascr$ is a graded sheaf-algebra on $X$
which resembles the symmetric algebra of a locally free sheaf of rank two on $X$.

In order to explain this definition we need a notion of locally free
sheaf in $\shbimod(X-X)$. We say that $\Escr\in \shbimod(X-X)$ is
locally free (of rank $n$) if $\pr_{1\ast}\Escr$ and
$\pr_{2\ast}\Escr$ are locally free (of rank $n$).
If $\Escr\in\shbimod(X-X)$ then we may define the tensor algebra
$T_X\Escr$ in the obvious way. If $\Escr$ is locally free of  rank two then in
\cite{Pat2,Pat1,VdB11}  a non-commutative symmetric algebra of rank two associated to
$\Escr$ is defined as a graded sheaf-algebra of the form
$T_X\Escr/(\Qscr)$ where $\Qscr\subset \Escr\otimes\Escr$ is
$\Qscr$ is locally free of rank one. While this is a reasonable
definition there are some problems with it.
\begin{itemize}
\item It is not so easy to find suitable $\Qscr$ inside $\Escr\otimes
  \Escr$ (see the complicated computations in \cite{VdB11}).
\item The dependence of $\QGr(T_X\Escr/(Q))$ on $\Qscr$ has not been
  made clear.  
\end{itemize}
In this paper we solve these problems by showing
that $\Qscr$ is actually superfluous (!) if $X$ is smooth. In other words
the theory can be
set up in a manner which does  not depend
  on an additional choice of  $\Qscr$. 

  We need the concept of a sheaf-$\ZZ$-algebra on $X$ . This is a
  sheaf algebra version of a usual $\ZZ$-algebra
  \cite{Bondal,Sierra}. Thus a sheaf $\ZZ$-algebra on~$X$ is
  defined by giving for $i,j\in \ZZ$ an object $\Ascr_{ij}$ in
  $\shbimod(X-X)$ together with ``multiplication maps''
  $\Ascr_{ij}\otimes \Ascr_{jk}\r\Ascr_{ik}$ and ``identity maps''
  $\Oscr_{X}\r \Ascr_{ii}$ satisfying the usual axioms. As in the
  graded case we may define abelian categories $\Gr(\Ascr)$ and
  $\QGr(\Ascr)$.

Let $\Escr$ be locally free of rank $n$. Then it is   easy to show that $-\otimes_{\Oscr_X} \Escr$ has
a right adjoint
$-\otimes_{\Oscr_X}\Escr^\ast$, where 
$\Escr^\ast\in \shbimod(X-X)$ is also locally free of
rank $n$  (this depends on $X$ being smooth).
Repeating this construction,  we may define $\Escr^{\ast 2}
=\Escr^{\ast\ast}$ by requiring that 
$-\otimes_{\Oscr_X}\Escr^{\ast\ast} $ is the right adjoint of 
$-\otimes_{\Oscr_X}\Escr^\ast$. By induction we define  $\Escr^{\ast 0}=\Escr$, $\Escr^{\ast m+1}
=(\Escr^{\ast m})^\ast$ for $m\ge
0$ and by  considering left adjoints we may define
$\Escr^{\ast m}$ for $m<0$. 

Standard properties of adjoint
functors yield  a bimodule inclusion
$
i_m:\Oscr_{X}\hookrightarrow \Escr^{\ast m}\otimes
\Escr^{\ast (m+1)},
$.

We now define $\SS(\Escr)$ as the $\ZZ$-algebra which satisfies
\begin{itemize}
\item[(a)]
 $\SS(\Escr)_{mm}=\Oscr_{X}$;
\item[(b)]
$\SS(\Escr)_{m,m+1}=\Escr^{\ast m}$;
\item[(c)] 
$\SS(\Escr)$ is freely generated by the $\SS(\Escr)_{m,m+1}$,
  subject to the relations given by the images of $i_m$.
\end{itemize}
\begin{definition}
\label{ref-1.1-0}
The non-commutative $\PP^1$-bundle $\PP(\Escr)$ on $X$ associated to
the locally free sheaf bimodule of rank two $\Escr$ on $X$ is the
category 
$\QGr(\SS(\Escr))$.
\end{definition}

It is easy to see that if $\Escr$ is an ordinary commutative
vector bundle of rank 2 on $X$ then
$\Gr(S_X(\Escr))\cong\Gr(\SS(\Escr))$. Thus the notion of a non-commutative
$\PP^1$-bundle if a generalization of the commutative one.
This is no
longer true in higher rank but even  then the algebra
$\SS(\Escr)$ could  be interesting in its own right. 

We will  show (see \S\ref{ref-4.2-25}) that if $\Escr\in\shbimod(X-X)$ is locally free of
rank $n$ and $\Qscr\subset
\Escr\otimes\Escr$ is of rank one and satisfies a
suitable non-degeneracy condition then
$\Gr(T_X\Escr/(\Qscr)=\Gr(\SS(\Escr))$. This shows that the current
definition of $\PP^1$-bundles is indeed a generalization of the
earlier one.

Let us now give a more detailed description of the content of this
paper. Our first main result is the following.
\begin{theorem} If $\Escr$ is  locally free of rank two then $\SS(\Escr)$ is a noetherian
sheaf-$\ZZ$-algebra in the sense that
$\Gr(\SS(\Escr))$ is a locally noetherian Grothendieck category.
\end{theorem}
To prove this we follow a standard approach (see \cite{ATV1})
which consists in defining a suitable quotient $\Dscr$ of  $\Ascr=\SS(\Escr)$
through the functor of point modules. The sheaf $\ZZ$-algebra $\Dscr$ will be noetherian by
construction and we will show that there is an invertible ideal
$\Jscr\subset \Ascr_{\ge 2}$ such that $\Dscr=\Ascr/\Jscr$. Then we
may conclude by invoking a suitable variant of the Hilbert basis
theorem.

Point modules over sheaf-($\ZZ$)-algebras have been defined in Adam Nyman's
PhD-thesis \cite{Nyman1} and he has shown that the corresponding functor is
representable (under suitable hypotheses). In particular it follows
from his results that the point functor of $\SS(\Escr)$ is representable
by $\PP_{X\times X}(\Escr)$. We reproduce the proof of this fact since
we need the exact nature of the bijections involved.

Our second main result is the following.
\begin{theorem}
\label{ref-1.3-1}
Assume that $Z$ is a Hirzebruch surface. Then every deformation of $Z$
is a non-commutative $\PP^1$-bundle over $\PP^1$.
\end{theorem}
For a precise definition of the notion of deformation we use we refer to
\S\ref{secdef} (which is based on \cite{VdBdef}).  The proof of Theorem \ref{ref-1.3-1} is
based on the observation that on $Z$ there are canonical exceptional line bundles
which may be lifted to any deformation. Imitating some
standard constructions in commutative algebraic geometry using the
resulting objects yields the desired result.

\medskip

After this paper was put on the arXiv the theory of non-commutative
$\PP^1$-bundles has been further developed. In \cite{Nyman5,Nyman3} it
was proved that they are $\Ext$-finite and satisfy a classical form of
Serre duality.  These papers use Theorem \ref{new--1} below. In return
the current proof of Theorem \ref{ref-1.3-1} uses some results from
\cite{Nyman5,Nyman3}.  

In~\cite{Mori1} it was shown that non-commutative
$\PP^1$-bundles share a number of geometric properties with their
commutative counterparts.  These results are stated in the language of
non-commutative algebraic geometry (where Grothendieck categories play
the role of spaces, see e.g.\ \cite{rosenberg,VdB19}).  In this
setting one may define a structure map $f:\PP(\Escr)\r X$ and Izuru
Mori shows that the fibers do not intersect. He also defines a certain
``quasi-section'' for $f$ and computes its self-intersection. In
\cite{Mori2} Izuru Mori computes the derived category of
non-commutative $\PP^1$-bundles.

In \cite{NymanChan} the authors attack the reverse question. They
generalize a standard characterization of ruled surfaces \cite{KM} to
the non-commutative case. Due to some new non-commutative phenomena that have
to be dealt with they do
not yet obtain a full analogue but nonetheless non-commutative
$\PP^1$-bundles appear as a basic example. Along the way the authors
prove that non-commutative $\PP^1$-bundles satisfy the Bondal-Kapranov
strengthening of Serre duality \cite{Bondal4} and are ``strongly
noetherian'' (which is important for the construction of Hilbert
schemes in this generality~\cite{AZ2}).
\section{Acknowledgement}
The author wishes to thank Daniel Chan, Colin Ingalls, Izuru Mori, Adman Nyman
and Paul Smith for useful discussions. In addition he thanks the anonymous
referee for his careful reading of the manuscript. 
\section{Notations and conventions}
Unless otherwise specified all schemes below will be of finite type
over a field $k$. 
\section{Sheaf-bimodules}
\label{ref-3-2}
\subsection{Generalities}
In the current and the next section we recapitulate the definition of
sheaf-bimodules from \cite{VdB11} and we give additional properties.
Since we will need to work with certain families of objects it will be
convenient to develop the material over a base-scheme $S$. In the
applications 
 we  will assume  $S=\Spec k$.

Below $S$ is a scheme and  $\alpha:X\r S$, $\beta:Y\r S$, $\gamma:Z\r
S$ will  be
$S$-schemes. An $S$-central coherent $X-Y$-sheaf-bimodule $\Escr$ is by definition a coherent
$\Oscr_{X\times_S Y}$-module such that the  support of $\Escr$ is finite
over
both $X$ and $Y$. We denote the corresponding abelian
category  by
$\shbimod_S(X-Y)$. More generally an  $S$-central
$X-Y$-sheaf-bimodule will be a quasi-coherent sheaf on $X\times_S Y$
which is a filtered direct limit of objects in
$\shbimod_S(X-Y)$. We denote the corresponding category by
$\ShBimod_S(X-Y)$. 
An object $\Escr$ in $\ShBimod_S(X-Y)$ 
defines a right exact functor $-\otimes_{\Oscr_X} \Escr :
\Qch(X)\r
\Qch(Y)$ commuting with direct sums via 
 $\pr_{2\ast}(\pr_1^\ast(-)\otimes_{\Oscr_{X\times_S
     Y}}\Escr)$. If $\Fscr$ is an object in 
 $\ShBimod_S(Y-Z)$
then the tensor product
$\Escr\otimes_{\Oscr_Y} \Fscr$ is defined as
$\pr_{13\ast}(\pr_{12}^\ast\Escr\otimes_{\Oscr_{X\times Y\times Z}}
\pr_{23}^\ast\Fscr)$. It is easy to show that this definition yields all
the expected properties (see \cite{VdB11}). 

Now  assume that we have finite $S$-maps
$u:W\r X$ and  $v: W\r Y$. If $\Uscr$ is a quasi-coherent
$\Oscr_W$-module then  we  denote the $X-Y$-bimodule 
$(u,v)_\ast\Uscr$ by ${}_u\Uscr_v$.  Any bimodule
$\Escr$ can be presented in this form
since we may take $W$ to be the scheme-theoretic support of $\Escr$. 
{}From the definition it is easy to check that $-\otimes
{}_u\Uscr_v=v_\ast(u^\ast(-)\otimes_{\Oscr_W} \Uscr)$.

It is useful to know that the functor $-\otimes_{\Oscr_X}\Escr$
actually determines $\Escr$. Let us  define $\Bimod(X-Y)$ as the category of right
exact functors $\Qch(X)\r \Qch(Y)$ commuting with direct
sums (this is equivalent to the definition in \cite{VdB19}).  
Then we have a functor
\[
F:\ShBimod_S(X-Y)\r \Bimod(X-Y)
\]
which sends $\Escr$ to the functor $-\otimes_{\Oscr_X} \Escr$. We have
the following result.
\begin{lemmas} 
\label{ref-3.1.1-3}
The functor $F$ is fully faithful.
\end{lemmas} 
\begin{proof}
We have to show how to reconstruct $\Escr$ from the functor
$-\otimes_{\Oscr_X} \Escr$.  

Choose an affine open covering
$X=\bigcup_i U_i$ and let $u_i:U_i\r X$, $u_{ij}:U_i\cap U_j \r X$ be the
inclusion maps. 

Assume that $H:\Qch(X)\r \Qch(Y)$ is a right exact functor commuting
with direct sums. Then $H(u_{i\ast} \Oscr_{U_i})$ will be a
quasi-coherent sheaf on $Y$ with an 
$\Oscr_X(U_i)$ structure. There is a corresponding quasi-coherent
sheaf $H_i$
on $U_i\times_S Y$.

In a similar way we find quasi-coherent sheaves $H_{ij}$ on $(U_i\cap
U_j) \times_S Y$ together with maps $H_i\mid U_i\cap U_j\r H_{ij}$.  We
define $\Fscr=\ker(\oplus_i u_{i\ast} H_i \r \oplus_{i\neq j}
u_{ij\ast}H_{ij})$. It is easy to see that if $H=-\otimes_{\Oscr_X}\Escr$ then
$\Fscr=\Escr$. 
\end{proof}
It would be interesting to give a more precise characterization of the
essential image of the functor $F$. One useful observation is that if $\Escr\in
\ShBimod_S(X-Y)$ then $-\otimes_{\Oscr_X} \Escr$ preserves exactness
of short exact sequence of vector bundles. This leads to the following
example.
\begin{examples}
Let $S=\Spec k$, $X=\PP^1$ and let $H:\Qch(X)\r \Qch(X)$ be the
functor given by $H(\Mscr)=\Oscr_{\PP^1} \otimes_k H^1(X,\Mscr)$. Then
$H$ does not preserve exactness of
\[
0\r \Oscr_{\PP^1}(-2)\r \Oscr_{\PP^1}(-1)^2 \r \Oscr_{\PP^1} \r 0
\]
and hence it is not in the essential image of $F$.

If we compute $\Fscr$ as in the proof of Lemma \ref{ref-3.1.1-3} then we
find $\Fscr=0$ which gives another reason why $H$ is not in the essential
image of $F$.
\end{examples}
A partial result in this context has been obtained by Nyman in \cite{Nyman4}.
\begin{definitions} An object $\Escr$ in $\shbimod_S(X-Y)$ is locally free on
  the left (right) (of rank $n$) if $\pr_{1\ast} \Escr$ ($\pr_{2\ast} \Escr$) is locally
  free on $X$ ($Y$) (of rank $n$).
\end{definitions}
 The following lemma shows that tensor products of
locally free bimodules behave as they should.
\begin{lemmas}
\label{rklemma}
Assume that $\Escr\in\shbimod_S(X-Y)$ and $\Fscr\in \shbimod_S(Y-Z)$
are locally free on the left. Then $\Escr\otimes_{\Oscr_Y} \Fscr$ is
also locally free on the left. Furthermore if $\Escr$ and $\Fscr$ have
constant rank on the left then so does $\Escr\otimes_{\Oscr_Y} \Fscr$
and the left rank of $\Escr\otimes_{\Oscr_Y} \Fscr$ is the product of the left ranks
of $\Escr$ and $\Fscr$.
\end{lemmas}
\begin{proof} 
As above we may assume $\Escr={}_u\Uscr_v$, $\Fscr={}_p \Vscr_q$ where $\Uscr$ is a coherent 
$W$-module for finite maps $u:W\r X$, $v:W\r Y$. Then $
\pr_{1\ast}(\Escr\otimes_{\Oscr_Y} \Fscr)=u_\ast(\Uscr\otimes_{\Oscr_W} v^\ast p_\ast \Vscr)$.
Here $\Vscr'=v^\ast p_\ast \Vscr$ is a locally free $\Oscr_W$-module. Thus we have
to show that if $u:W\r X$ is a finite map and $\Uscr,\Vscr'$ are coherent $\Oscr_W$-modules
such that $\Vscr'$ is locally free and $u_\ast \Uscr$ is locally free then $u_\ast(\Uscr\otimes_{\Oscr_W} \Vscr')$ is locally free.  Since the question is local on $X$ we may reduce
to the case that $X$ is affine. Then $W$ is affine as well and hence 
$\Vscr'$ is a direct summand of a free $\Oscr_W$-module. So we reduce to the case 
$\Vscr'=\Oscr_W$ which is obvious.


It is sufficient to prove the assertion on the rank for all pullbacks $\Spec l\r S$
for $l$ algebraically closed. Hence we may assume that $S=\Spec l$ with
$k$ algebraically closed.

Now let $m$, $n$ be respectively the left rank of $\Escr$ and $\Fscr$.
We have to show that $
  \length(\Oscr_x\otimes_{\Oscr_X} \Escr\otimes_{\Oscr_Y}\Fscr)=mn$ for
  all closed points $x\in X$. Since $\Oscr_x\otimes_{\Oscr_X}\Escr$ is an
  extension of $m$ objects of the form $\Oscr_{y_i}$ for some $y_i\in
  Y$, this is clear.
\end{proof}
In the sequel we will use the following lemma to show that certain
sheaves are locally free.
\begin{lemmas} 
\label{ref-3.1.5-4}
Assume that $\psi:R\r S$ is a local ring homomorphism between noetherian commutative local
  rings with maximal ideals $m,n$. Let $u:M\r N$  be a morphism
  between finitely generated $S$ modules where $N$ is in addition flat
  over 
   $R$. Assume that $u\otimes_{R} R/m$ 
  is injective with $S/mS$-free cokernel. Then $u$ is also injective
  with 
  $S$-free cokernel.
\end{lemmas}
\begin{proof}
Let $C$ be the cokernel of $u$.  By
hypotheses $C/m C$ is free over $S/mS$. Choose an
isomorphism $(S/mS)^k\r C/m C$ and lift this to a map
$\theta:S^k\r C$.  Let $T$ its cokernel. Tensoring with $R/m$ yields $T/mT=0$. Since
$\psi(m)\subset n$ we obtain $T=0$ by Nakayama's lemma. Now factor
$\theta$ through a map $\theta':S^k\r N$ and let $K$ be the pullback
of $\theta'$ and $u$. Thus we have an exact sequence:
\[
0\r K\r M\oplus S^k\xrightarrow{(u,\theta')} N \r 0
\]
Since $N$ is flat over $R$ this sequence remains exact if we tensor
with $R/mR$. Since $(S/mS)^k$ is isomorphic to coker $u\otimes_R
R/m$  we deduce that $K/mK=0$. By Nakayama we obtain $K=0$. This
clearly implies what we want. 
\end{proof}
If $\alpha$ is smooth then we will say that $\alpha$ is
\emph{equidimensional} if the fibers of $\alpha$ are equidimensional
and if furthermore they all have the same dimension. We will say that
$\alpha$ is of \emph{relative dimension $n$} if it is equidimensional and if
all fibers have dimension $n$.

The following result will be very convenient:
\begin{propositions}
\label{ref-3.1.6-5}
Assume that $\alpha,\beta$ are smooth and equidimensional of the same
relative dimension. Then $\Escr\in \shbimod(X-Y)$ is locally free on the
left if and only it is locally free on the right.
\end{propositions}
\begin{proof}
Assume that $\Escr$ is locally free on the left. We will show that it
is also locally free on the right. 
First consider the case that $S=\Spec k$. Then $X$ and
$Y$ are regular of the same dimension. As above we may assume that
$\Escr={}_\delta \Uscr_\epsilon$ for finite maps $\delta:W\r X$,
$\epsilon:W\r Y$. 
We then have the following chain of implications:
\begin{align*}
\text{$\delta_{\ast} \Uscr$ is locally free} &\Rightarrow
\text{$\delta_{\ast} \Uscr$ is maximal Cohen-Macaulay}\\
 &\Rightarrow \text{$\Uscr$ is maximal Cohen-Macaulay on $W$}\\
&\Rightarrow\text{$\epsilon_{\ast} \Uscr$ is maximal Cohen-Macaulay}\\
&\Rightarrow \text{$\epsilon_{\ast} \Uscr$ is locally free}
\end{align*}
The last implication follows from the fact that $Y$ is regular. 

Now consider the case where $S$ is general. From the hypotheses that $\delta_\ast \Uscr$
is locally free over $X$ we obtain that $\Uscr$ is flat over $S$ and hence $\epsilon_\ast\Uscr$ is also flat over $S$ (since $\epsilon$ is finite). 

Thus $\epsilon_{\ast} \Uscr$ is flat over $S$. Since $\epsilon$ is
finite the formation of $\epsilon_{\ast} \Uscr$ commutes with base
change.  By the above discussion we know that for every $s\in S$ we
have that $\epsilon_\ast (\Uscr_s)$ is locally free over $Y_s$.  Then
lemma \ref{ref-3.1.5-4} with $M=0$ shows that $\epsilon_{\ast} \Uscr$
itself is locally free.
\end{proof}
Below we assume that $\alpha:X\r S$, $\beta:Y\r S$, $\gamma:Z\r S$
are smooth and equidimensional of the same relative dimension.

Now assume that $\Escr$ is an  object in $\shbimod_S(X-Y)$ which
is locally free on the left (and hence on the right). We will
define/construct the right and left duals $\Escr^\ast$, ${}^\ast\Escr$ to
$\Escr$. For brevity we restrict the discussion below to the right dual. Everything
has obvious analogues for the left dual. 

We want $\Escr^\ast\in\shbimod_S(Y-X)$ and in addition we should have
\[
\Hom_Y(\Ascr\otimes_{\Oscr_X} \Escr,\Bscr)= \Hom_X(\Ascr,
\Bscr\otimes_{\Oscr_Y} \Escr^\ast)
\]
According to lemma \ref{ref-3.1.1-3} this property defines $\Escr^\ast$ up to unique
isomorphism, if it exists.

We now describe $-\otimes_{\Oscr_Y} \Escr^\ast$ more
precisely. With the same notations as before we assume
$\Escr={}_{u}\Uscr_{v}$ where $\Uscr\in\coh(W)$. Let us
denote with $v^!$ the right adjoint to~$v_\ast$. Then it is
easy to verify that  one has
\[
\Escr^\ast={}_v\HHom_W(\Uscr,v^!\Oscr_Y)_u
\]
from which in particular we deduce 
\begin{equation}
\label{ref-3.1-6}
\pr_{1\ast}(\Escr^\ast)\cong \pr_{2\ast}(\Escr)^\ast
\end{equation}
Thus the left structure of $\Escr^\ast$ is given by the dual of the
right structure of $\Escr$.

Let $Rv^!$ be the right derived functor to $v^!$ (note that
this is somewhat at variance with the usual definitions). Then it is
clear that we also have
\begin{equation}
\label{ref-3.2-7}
\Escr^\ast={}_v\HHom_W(\Uscr,Rv^!\Oscr_Y)_u={}_v\uRHom_W(\Uscr,Rv^!\Oscr)_u
\end{equation}
Furthermore if $\omega_{X/S}$ denotes the relative dualizing complex
then we have $Rv^!(\Oscr_Y)=\omega_{W/S}\otimes_{\Oscr_W}
v^\ast \omega_{Y/S}^{-1}$ from which we deduce
\begin{equation}
\label{ref-3.3-8}
\Escr^\ast= \omega_{Y/S}^{-1} \otimes_{\Oscr_Y}{}_v(\Uscr^D)_u
\end{equation}
where $(-)^D$ denotes the Cohen-Macaulay dual. By symmetry we have a
similar formula
\begin{equation}
\label{ref-3.4-9}
{}^\ast\Escr=  {}_v(\Uscr^D)_u\otimes_{\Oscr_X}\omega_{X/S}^{-1}
\end{equation}
where ${}^\ast\Escr$ is defined as $\Escr^\ast$ but using left
adjoints. 
\begin{lemmas} 
\label{ref-3.1.8-10}
We have $\Escr^{\ast\ast} =\omega^{-1}_{X/S}
  \otimes_{\Oscr_X} \Escr  \otimes_{\Oscr_Y} \omega_{Y/S}$
\end{lemmas}
\begin{proof} The author learned this beautiful formula from  notes by 
  Kontsevich \cite{Ko} where it is shown that it  holds more generally in the setting of
  derived categories. In our current setting it follows trivially from
  \eqref{ref-3.3-8}. 
\def\qed{}\end{proof}
\begin{corollarys}
\label{ref-3.1.9-11}
 The left rank of $\Escr$ equals the right rank of
  $\Escr^\ast$ and vice-versa.
\end{corollarys} 
\begin{proof} 
According to \eqref{ref-3.1-6}  the left
structure of $\Escr^\ast$ is given by the ordinary vector bundle dual
of the right structure of $\Escr$. Thus the right rank of $\Escr$
equals the left rank of $\Escr^\ast$. In the same way we find that the
right rank of $\Escr^\ast$ equals the left rank  of
$\Escr^{\ast\ast}$. Now from lemma \ref{ref-3.1.8-10} we easily obtain
that the left rank of $\Escr^{\ast\ast}$ equals the left rank of
$\Escr$ which finishes the proof.
\end{proof}
The following lemma will be used many times.
\begin{lemmas} \label{ref-3.1.10-12} The formation of $(-)^\ast$ is compatible with
  base change for locally free coherent sheaf-bimodules. 
\end{lemmas}
\begin{proof}
If $\Escr$ is a locally free coherent sheaf-bimodule on
$X$ and we have a base extension $T\r S$ then using the formula \eqref{ref-3.3-8} we see that there
is at least a map of sheaf-bimodules $
(\Escr^\ast)_T \r (\Escr_T)^\ast$. Then by looking at the left or right structure we see that this
map is an isomorphism. 
\end{proof}

Using standard properties of adjoint functors together with lemma
\ref{ref-3.1.1-3} we obtain canonical maps in $\ShBimod_S(X-Y)$ 
\begin{gather*}
i:\Oscr_X\r \Escr\otimes_{\Oscr_Y} \Escr^\ast\\
j:\Escr^\ast\otimes_{\Oscr_X} \Escr\r \Oscr_Y 
\end{gather*}
In the sequel we will need some properties of these maps. 
\begin{propositions}
\begin{enumerate} \item
$i$ is injective and its cokernel is locally free.
\item $j$ is surjective (and hence its kernel is trivially locally
  free).
\end{enumerate}
\end{propositions}
\begin{proof}
We only consider (1). (2) is similar. With a similar method as the one
that was used in the proof of Proposition \ref{ref-3.1.6-5} it suffices to
prove this in the case that $S=\Spec k$. If we
restrict to this case then it is sufficient to prove that for all
closed points $x\in X$ the map 
\[
\Oscr_x\r \Oscr_x\otimes_{\Oscr_X} \Escr\otimes_{\Oscr_Y} \Escr^\ast
\]
is non-zero. Now this map is obtained by adjointness from the identity
map
\[
\Oscr_x\otimes_{\Oscr_X} \Escr
\r 
\Oscr_x\otimes_{\Oscr_X} \Escr
\]
Since this map is obviously non-zero we are done. 
\end{proof}
Below it will be convenient to have a slight generalization of the
relationship that exists between members of a pair
$(\Escr,\Escr^\ast)$.

 Therefore we make the following definition.
\begin{definitions}
$\Qscr\in\shbimod_S(X-Z)$ is invertible if there exists $\Qscr^{-1}\in
\shbimod_S(Z-X)$ together with isomorphisms
$\Qscr\otimes_{\Oscr_Z}\Qscr^{-1}\cong \Oscr_X$ and
$\Qscr^{-1}\otimes_{\Oscr_X} \Qscr\cong\Oscr_Z$.
\end{definitions}
Using the results in \cite{AVdB} or \cite{AZ} one obtains that 
$\Qscr\in \shbimod_S(X-Z)$ is invertible if and only if 
$\Qscr\cong {}_{\id_X}(\Lscr)_\beta$ where $\Lscr\in\Pic(X)$ and $\beta$
is an $S$ isomorphism between $X$ and $Z$.
\begin{definitions}  Let $\Escr,\Fscr$ be  locally
  free objects respectively in $\shbimod(X-Y)$ and $\shbimod(Y-Z)$.
  Assume that $\Qscr$ is an invertible object in
  $\shbimod(X-Z)$ and assume furthermore that $\Qscr$ is contained in
  $\Escr\otimes_{\Oscr_Y} \Fscr$. We say that $\Qscr$ is
  non-degenerate if the following composition
\[
\Escr^\ast\otimes_{\Oscr_X} \Qscr\r \Escr^\ast\otimes_{\Oscr_X}
\Escr\otimes_{\Oscr_Y} \Fscr\r \Fscr
\]
is an isomorphism.
\end{definitions}
Clearly if $\Qscr$ is non-degenerate in $\Escr\otimes_{\Oscr_Y} \Fscr$
then we have 
\begin{equation}
\label{ref-3.5-13}
\Escr^\ast\cong \Fscr\otimes_{\Oscr_Z}\Qscr^{-1}
\end{equation}

\subsection{Sheaf algebras and sheaf $\ZZ$-algebras}
\label{ref-3.2-14}
In this section the notations will be as in the previous section. It
is clear that $\ShBimod_S(X-X)$ is a monoidal category so we can
routinely define algebras and $I$-algebras in this category (see
\cite{Bondal} for the definition of ordinary $\ZZ$-algebras. If we
replace the indexing set $\ZZ$ by an arbitrary set $I$ then we obtain
the notion of an $I$-algebra). We will call these ($S$-central)
sheaf-algebras and ($S$-central) sheaf-$I$-algebras.  For example a
sheaf-algebra on $X$ is an object $\Ascr$ in $\ShBimod_S(X-X)$
together with a multiplication map
$\Ascr\otimes_{\Oscr_X}\Ascr\r\Ascr$ and a unit map $\Oscr_X\r \Ascr$
having the usual properties. If $\Ascr$ is a sheaf-algebra on $X$ then
we define $\Mod(\Ascr)$ as the category consisting of objects in
$\Qch(X)$ together with a multiplication map $\Mscr\otimes_{\Oscr_X}
\Ascr\r \Mscr$, again satisfying the usual properties. In the same way
we may define $\ShBimod(\Ascr-\Ascr)$. This and similar notions will
be used routinely in the sequel. We leave the obvious definitions to
the reader.

The previous paragraph makes clear what we mean by a
 sheaf-$I$-algebra on $X$. However in the sequel we will use this
 notion in somewhat greater generality. So we will discuss this next.

 Assume that $\Xi$ is a
family of $S$ schemes $\alpha_i:X_i\r S$ indexed by $i\in I$.  A
sheaf $I$-algebra on  $\Xi$
 is defined by giving for $i,j\in I$ an object
$\Ascr_{ij}$ in $\ShBimod_S(X_i-X_j)$ together with ``multiplication
maps'' $\Ascr_{ij}\otimes_{\Oscr_{X_j}} \Ascr_{jk}\r\Ascr_{ik}$ and an ``identity
map'' $\Oscr_{X_i}\r \Ascr_{ii}$ satisfying the usual 
axioms. 

If $\Ascr$ is a sheaf-$\Xi$-algebra then an \emph{$\Ascr$}-module is a
formal direct sum  $\oplus_{i\in I}\Mscr_i$ where $\Mscr_i\in
\Qch({X_i})$ together with multiplication maps $\Mscr_{i}\otimes_{\Oscr_{X_i}}
\Ascr_{ij}\r \Mscr_j$, again satisfying the usual axioms. We denote the
category of $\Ascr$-modules by $\Gr(\Ascr)$. It is easy to see that
$\Gr(\Ascr)$ is a Grothendieck category. 

Unless otherwise specified we
will now assume that $I=\ZZ$ even though some (but not all) notions
below make sense more generally.
We will say that $\Ascr$ is noetherian if $\Gr(\Ascr)$ is a locally
noetherian abelian category. In the case that
$\Ascr$ is noetherian we borrow a number of definitions from \cite{AZ}.
Let $M\in \Gr(\Ascr)$. We say that $M$ is is \emph{left}, resp.\ 
\emph{right bounded} if $M_i=0$ for $i\ll 0$ resp. $i\gg 0$. We say that
$M$ is \emph{bounded} if $M$ is both left and right bounded. We say
$M$ is \emph{torsion} if it is a direct limit of right bounded objects.
We denote the corresponding category by $\Tors(\Ascr)$. Following
\cite{AZ} we also put $\QGr(\Ascr)=\Gr(\Ascr)/\Tors(\Ascr)$. Furthermore we define
the following functors. $\tau:\Gr(\Ascr)\r \Tors(\Ascr)$ is the torsion
functor associated to $\Tors(\Ascr)$; $\pi :\Gr(\Ascr)\r \QGr(\Ascr)$ is the
quotient functor; $\omega:\QGr(A)\r \Gr(A)$ is the right adjoint to
$\pi$ and finally $\widetilde{(-)}=\omega\pi$.

In these notes we will use the convention that if $\text{Xyz}$ is an abelian
category then $\text{xyz}$ denotes the full subcategory of $\text{Xyz}$ whose objects
are given by the noetherian objects. Following this convention we
introduce $\qgr(\Ascr)$ and $\tors(\Ascr)$. Note that if $M\in \tors(\Ascr)$
then $M$ is right bounded, just as in the ordinary graded case. 
It is also easy to see that
$\qgr(\Ascr)$ is equal to $\gr(\Ascr)/\tors(\Ascr)$. 
We put $\Ascr_{\ge l}=\oplus_{j- i\ge l} \Ascr_{ij}$ and similarly $\Ascr_{\le
  l}=\oplus_{j-i\le l} A_{ij}$. $\Ascr_{\ge 0}$ and $\Ascr_{\le 0}$
are both sheaf-$\ZZ$-subalgebras of $\Ascr$ and $\Ascr_{\ge l}$ and
$\Ascr_{\le l}$ are sheaf-bimodules over respectively $\Ascr_{\ge 0}$ and $\Ascr_{\le 0}$.

We say that $\Ascr$ is positive  if $\Ascr=\Ascr_{\ge 0}$.
\begin{lemmas}\cite{NVO} $\Ascr$ is noetherian if and only if $\Ascr_{\ge 0}$
  and $\Ascr_{\le 0}$ are noetherian.
\end{lemmas}
We will use the following generalization of the Hilbert
basis-theorem. 
\begin{lemmas}
\label{ref-3.2.2-15}  Assume that $\Ascr$ is positive and let $I\subset
  \Ascr_{\ge 1}$  be an invertible ideal in $\Ascr$ (that is an
  invertible object in $\ShBimod(\Ascr-\Ascr)$ which is contained 
  in $\Ascr$). If $\Ascr/I$ is noetherian then so is $\Ascr$.
\end{lemmas}
$\Ascr$ is said to be strongly graded if the canonical map
$\Ascr_{ij}\otimes_{\Oscr_{X_j}} \Ascr_{jk}\r Ascr_{ik}$ is surjective for all
$i,j,k$. We have \cite{NVO}
\begin{lemmas} If $\Ascr$ is strongly graded then the restriction
  functor $\Gr(\Ascr)\r \Mod(\Ascr_{ii}): M\mapsto M_i$ is an
  equivalence of categories for all $i$.
\end{lemmas}

An interesting fact about sheaf-$\ZZ$-algebras is that they admit a useful
form of twisting. 
 Let $\Ascr$ be a sheaf-$\ZZ$-algebra over $\Xi$ and let
 $\Xi'=(X'_i)_{i\in \ZZ}$ be another family of $S$-schemes.
Let $\Tscr_{i}$ be invertible objects in 
$\ShBimod_S(X_i-X'_i)$. Define the sheaf-$\ZZ$-algebra
$\Bscr$ via
\[
\Bscr_{ij}=\Tscr_i^{-1}\otimes_{\Oscr_{X_i}} \Ascr_{ij}
\otimes_{\Oscr_{X_j}} \Tscr_j
\]
It is easy to see that the functor
\[
\oplus_i \Mscr_i \mapsto \oplus_i \Mscr_i
 \otimes_{\Oscr_{X_i}}\Tscr_i \]
defines an equivalence $\Gr(\Ascr)\cong \Gr(\Bscr)$.

\subsection{Ampleness}
If $\Xi=(\alpha_i:X_i\r S)_{i\in \ZZ}$ and $\Omega=(\beta_i:Y_i\r
S)_{i\in \ZZ}$ are collections of $S$-schemes then a map
$\gamma:\Omega\r \Xi$ is a collection of maps $(\gamma_i:Y_i\r
X_i)_{i\in \ZZ}$ such that
$\alpha_i\gamma_i=\beta_i$. Assume now  that the following
condition holds for $\gamma$:
\begin{itemize}
\item[(C)] Let $i,j\in\ZZ$ be arbitrary and let $Z$ be an arbitrary closed
  subset of $Y_i\times_S Y_j$ which is finite over both factors. Then
  the image of $Z$ in $X_i\times_S X_j$ is also finite over both
  factors. 
\end{itemize}
\begin{examples} Here is an example why this condition is not vacuous
  even if $Y_i\r X_i$ is proper. Let $S=\Spec k$ and let $(E,+)$ be an
  elliptic curve over $k$. Assume $Y_i=Y_j=E\times E$ and $X_i=X_j=E$
  where $\gamma_i$ is the projection on the first factor in $E\times E$.  Let $Z\subset
  (E\times E)\times (E\times E)$ be the graph of the automorphism
  $E\times E\r E\times E:(x,y)\mapsto (x+y,y)$. Then the projection of
  $Z$ on $E\times E$ is $E\times E$ and hence is not finite over both
  factors.
\end{examples}
If $\Bscr$ is a  sheaf-$\ZZ$-algebra on $\Omega$ and $\gamma$ satisfies (C) then we may define
sheaf $\ZZ$-algebra $\gamma_\ast(\Bscr)$ on $\Xi$ by 
\[
\gamma_\ast(\Bscr)_{ij}=(\gamma_i,\gamma_j)_\ast(\Bscr_{ij})
\]
There is a canonical
functor $\gamma_\ast:\Gr(\Bscr)\r \Gr(\gamma_\ast\Bscr):
\oplus_i\Mscr_i\mapsto \oplus_i \gamma_{i,\ast}\Mscr_i$. This functor factors
through a functor $\bar{\gamma}_\ast:\QGr(\Bscr)\r \QGr(\gamma_\ast\Bscr)$.
In the sequel
we will study the properties of this functor in some special cases. 

Let us now assume that $\Bscr$ is a positive sheaf-$\ZZ$-algebra on $\Omega$ such that
all 
$\Bscr_{ij}$ are coherent. Assume furthermore that all $\gamma_i$ are proper.
Examining 
\cite{AVdB,VdB11} leads to the following notion.
\begin{definitions} 
\label{ref-3.3.1-16}
$\Bscr$ is  ample for $\gamma$ if the following
  conditions hold
\begin{enumerate}
\item $\Bscr$ is noetherian.
\item For every $i\in \ZZ$ and $\Mscr\in\coh(Y_i)$ we have that
  $\Mscr\otimes_{\Oscr_{Y_i}}\Bscr_{ij}$ is relatively generated by
  global sections for the map $\gamma_j$ for $j\gg 0$. 
\item For every $i\in \ZZ$, $k>0$ and $\Mscr\in\coh(Y_i)$ we have that
  $R^k\gamma_{j\ast}(\Mscr\otimes_{\Oscr_{Y_i}}\Bscr_{ij})=0$ for
  $j\gg 0$.
\end{enumerate}
\end{definitions}
Generalizing \cite{AZ,AVdB,VdB11} we then obtain:
\begin{theorems}
\label{ref-3.3.2-17}
Assume that condition (C)
holds and that all $\gamma_i$ are proper. Assume furthermore that $\Bscr$ is ample for $\gamma$.  Then $\bar{\gamma}_\ast$ is
an equivalence of categories. In addition $\gamma_\ast(\Bscr)$ is
noetherian and the functor $\gamma_\ast$ preserves noetherian
objects. 
\end{theorems}

\subsection{Point modules}
Point modules over sheaf-($\ZZ$-)algebras have been introduced by Adam
Nyman in his PhD-thesis \cite{Nyman1}. We reproduce his definition
below.

We first introduce another notion of local freeness. If $\alpha:X\r S$ is an $S$-scheme and
$P\in \coh(X)$ then we say that $P$ is coherent over $S$ if the
support of $P$ is finite over $S$.

We say that $P$ is locally free (of rank
$n$) over $S$ if  $P$ is coherent over $S$ and $\alpha_\ast
P$ is locally free (of rank $n$). If $P$ is locally
free of rank one over $S$ then it is of the form $\zeta_\ast Q$ for a unique
section $\zeta:S\r X$ of $\alpha$ and $Q$ a line bundle on $S$. Using a
slight abuse of notation we write $P^{-1}$ for
$\zeta_\ast(Q^{-1})$. If $\alpha:X\r S$ and $\beta:Y\r S$
are $S$-schemes and  if $P_1\in\coh(X)$, $P_2\in\coh(Y)$ are locally
free of rank one over $S$ then so is 
\[
P_1\boxtimes_S P_2=\pr_1^\ast(P_1)\otimes_{\Oscr_{X\times_S Y}}
\pr_2^\ast(P_2)
\]
Note that if $P_1=\zeta_{1\ast}(Q_1)$ and $P_2=\zeta_{2\ast}(Q_2)$ then 
\[
P_1\boxtimes_S P_2=(\zeta_1,\zeta_2)_\ast(Q_1\otimes_{\Oscr_S} Q_2)
\]
We will need the following result.
\begin{lemmas}
\label{ref-3.4.1-18}
Assume that $\alpha:X\r S$ and $\beta:Y\r S$
are $S$-schemes
and let
$\Escr\in\ShBimod_S(X-Y)$. Let $P_0\in \coh(X)$, $P_1\in \coh(Y)$ be
locally free of
rank one over $S$. Then we have  canonical isomorphisms:
\begin{equation}
\label{ref-3.6-19}
\Hom_{\Oscr_Y}(P_0\otimes_{\Oscr_X} \Escr,P_1)
\cong\Hom_{\Oscr_{X\times_S Y}}(\Escr,P_0^{-1}\boxtimes_S P_1)
\end{equation}
Furthermore under this isomorphism, epimorphisms correspond to each
other.
\end{lemmas}
\begin{proof}
This is a direct computation. Let $P_0=\zeta_{0\ast}(Q_0)$,
$P_1=\zeta_{1\ast}(Q_1)$ where $\zeta_{1}:S\r X$, $\zeta_2:S\r Y$  are sections of
$\alpha$ and $\beta$ respectively. We have
\[
P_0\otimes_{\Oscr_X}
\Escr=\pr_{2\ast}(\pr_1^\ast\zeta_{0\ast}Q_0\otimes_{\Oscr_{X\times_SY }} \Escr)
\]
Thus we have
\[
\Hom_{\Oscr_Y}(P_0\otimes_{\Oscr_X} \Escr,P_1)
=
\Hom_{\Oscr_S}(\zeta_1^\ast\pr_{2\ast}(\pr_1^\ast\zeta_{0\ast}Q_0\otimes_{\Oscr_{X\times_S Y
    }} \Escr), Q_1)
\]
If we look at the following pullback diagram:
\[
\begin{CD}
X @>(\Id_X,\zeta_1\alpha) >> X\times_S Y\\
@V\alpha VV @V\pr_2 VV\\
S @>\zeta_1 >> Y
\end{CD}
\]
then we find
\begin{align*}
\zeta_1^\ast\pr_{2\ast}(\pr_1^\ast\zeta_{0\ast}Q_0\otimes_{\Oscr_{X\times_S Y
    }} \Escr)&=
\alpha_\ast (\Id_X,\zeta_1\alpha)^\ast(\pr_1^\ast\zeta_{0\ast}
    Q_0\otimes_{\Oscr_{X\times_S Y}} \Escr)\\
&=\alpha_\ast(\zeta_{0\ast}Q_0\otimes_{\Oscr_X}
    (\Id_X,\zeta_1\alpha)^\ast \Escr)\\
&=\alpha_\ast\zeta_{0\ast}(Q_0\otimes_{\Oscr_S}
    \zeta^\ast_0(\Id_X,\zeta_1\alpha)^\ast \Escr)\\
&=Q_0\otimes_{\Oscr_S} (\zeta_0,\zeta_1)^\ast(\Escr)
\end{align*}
We now compute
\begin{align*}
\Hom_{\Oscr_S}(Q_0\otimes_{\Oscr_S} (\zeta_0,\zeta_1)^\ast
\Escr,Q_1)&=
\Hom_{\Oscr_S}((\zeta_0,\zeta_1)^\ast
\Escr,Q_0^{-1}\otimes_{\Oscr_S} Q_1)\\
&=\Hom_{\Oscr_{X\times_S Y}}(
\Escr,(\zeta_0,\zeta_1)_\ast(Q_0^{-1}\otimes_{\Oscr_S} Q_1))\\
\\
&=\Hom_{\Oscr_{X\times_SY}}(
\Escr,P_0^{-1}\boxtimes_S P_1)\\
\end{align*}
To prove the claim about preservation of epimorphisms one simply
checks that epimorphisms are preserved in each individual
step. 
\end{proof}
Assume now that $\Ascr$ is a positively graded
sheaf-$\ZZ$-algebra on $\Xi$. Just as in the case of ordinary algebras
one may define a concept of point modules in $\Gr(\Ascr)$. 
\begin{definitions} 
An $m$-shifted  point module over $\Ascr$ is an $\Ascr$-module ${{P}}$
generated in degree $m$ such that 
for $n\ge m$ we have that $ {{P}}_n$ is locally free of rank one over $S$. A $0$-shifted point module will be
simply called a point-module. An extended point module over $\Ascr$ is
an $\Ascr$-module $P$ such that for all $m$, $P_{\ge m}$ is an
$m$-shifted point module.
\end{definitions}
To study point modules it will be convenient to introduce the notion
of  a truncated point module. Let $[m:n]=\{m,m+1,\ldots,n\}$ and let 
 $\Ascr_{[m:n]}=\oplus_{m\le i,j\le n} \Ascr_{ij}$. Clearly
$\Ascr_{[m:n]}$ is a $[m:n]$-algebra. There are obvious
restriction functors $\Gr(\Ascr)\r \Gr(\Ascr_{[m:n]})$ and
$\Gr(\Ascr_{[m:n]})\r \Gr(\Ascr_{[m':n']})$ when $m'\ge m, n'\le n$.

 We
define a $[m:n]$-truncated $\Ascr$ point module ${{P}}$  as an $\Ascr_{[m:n]}$-module
generated in degree $m$ such that 
for $n\ge i\ge m$ we have that ${{P}}_i$ is locally free of rank one.

It is natural to declare two (truncated, extended, shifted) point
modules $P,Q$ to be equivalent if  there exists a line bundle $\Lscr$ on $S$ such
that $Q_n=\alpha_n^\ast\Lscr\otimes_{\Oscr_{X_n}} P_n$

The main feature of (extended) point modules is that they define
certain sheaf-$\ZZ$-algebras 
which may be used to study $\Ascr$. Let $P$ be an extended
point module over $\Ascr$. Thus for every $i$ we have that $P_i$ is locally free of rank one
over $S$ and hence $P_i=\zeta_{i\ast}(Q_i)$ where $\zeta_i$ is a section of
$\alpha_i$ and $Q_i\in\Pic(S)$.

 We define $\Bscr_{ij}(P)=Q_i^{-1}\otimes_S Q_j$. 
Thus $\Bscr(P)=\oplus_{ij}\Bscr_{ij}(P)$ is a strongly graded sheaf-$\ZZ$-algebra on
$S$.  Let $\Omega=(S)_{i\in\ZZ}$ be the trivial constant system of $S$-schemes
and let $\zeta:\Omega\r \Xi$ be defined by $(\zeta_i)_i$. Then
the right $\Ascr$-module structure of $P$ yields us through
lemma \ref{ref-3.4.1-18} a surjective map $\Ascr_{m,n}\r \zeta_\ast\Bscr_{m,n}(P)$ and a
straightforward verification shows that this map is compatible with
multiplication. Hence we obtain a surjective map of sheaf-$\ZZ$-algebras
$\Ascr\r \zeta_\ast\Bscr(P)$.

In the sequel we will need families of the concepts that were
introduced above. If $\theta:W\r S$ is an $S$-scheme then we can consider
the base extended algebra $\Ascr_W$ which is just $\oplus_{m,n}
(\theta,\theta)^\ast(\Ascr_{m,n})$ where we have denoted the base
extension of $\theta$ to a map $X_{n,W}\r X_n$ also by $\theta$. We
define a family of point modules over $\Ascr$ parametrized by $W$ to
be a point module on $\Ascr_W$. Families of extended and truncated
point modules are defined in a similar way. 

Assume that $P$ is a family of extended point modules parametrized by $W$. Then
$\Bscr(P)$ is a $W$-central sheaf-$\ZZ$-algebra on $W$. As above we
have $P_i=\zeta_{i\at}(Q_i)$ where $Q_i\in \Pic(W)$ and $\zeta_i$ is a
section of $X_i\r X_{i,W}$. We may write $\zeta_i$ as $(\mu_i,\Id_W)$
with $\mu_i$ a map $W\r X_i$.
\begin{lemmas}
The image of $(\mu_i,\mu_{j})$ lies inside the support of $\Ascr_{ij}$.
\end{lemmas}
\begin{proof}
By the definition of a point module we have a surjective map
\[
P_i\otimes_{\Oscr_{X_{i,W}}}\Ascr_{W,ij}\r P_j
\]
which according to lemma \ref{ref-3.4.1-18} corresponds to a surjective map
\[
(\theta,\theta)^\ast(\Ascr_{ij}) \r P_i^{-1}\boxtimes_W P_j
\]
Thus the image of $(\zeta_i,\zeta_j)$ lies inside
$(\theta,\theta)^{-1} (\Supp \Ascr_{ij})$. It follows that the image of
$(\mu_i,\mu_j)=(\theta\circ \zeta_i,\theta\circ\zeta_j)$ lies inside $\Supp
\Ascr_{ij}$.  This proves what we want.
\end{proof}
\begin{corollarys}
\label{ref-3.4.4-20}
 Assume that $\theta:W\r S$ is proper and that all
  $\Ascr_{ij}$ are coherent. Then the $\mu_i$ are proper. Let
  $\Omega=(W)_{i\in\ZZ}$ be the constant system associated to $W$
  and let $\mu:\Omega\r \Xi$ be given by $(\mu_i)_i$. Then $\mu$
  satisfies $(C)$ and the map $\Ascr_W\r \Bscr(P)$ gives by
  adjointness rise to a map $\Ascr\r \mu_\ast \Bscr(P)$.
\end{corollarys}
\begin{proof} The map $\mu_i$ is the composition
  $W\xrightarrow{\zeta_i} X_{i,W}\xrightarrow{\theta} X_i$. The first map is a
  section and so it is a closed
  immersion. In particular it is proper. The second map is
  also proper since it is the base extension of a proper map. Thus
  $\mu_i$ is also proper. 
  
  Now we can verify (C). Since $(\mu_i,\mu_j)$ is proper  it is
  sufficient to verify that the image of 
  $(\mu_i,\mu_j)$ is finite on the left and right. This is clear since
  by the previous lemma this image is contained in the support of
  $\Ascr_{ij}$ and $\Ascr_{ij}$ was coherent by hypotheses.
\end{proof}

Equivalences among
families of point modules are defined in the same way as for ordinary
point modules (see above). For use in the sequel we introduce the following (somewhat
adhoc) notations.
\[
\begin{array}{ll}
\Points_{m,n}(W) & \text{equivalence classes of $[m:n]$-truncated point modules parametrized by
  $W$}\\
\Points_m(W) &\text{equivalence classes of $m$-shifted point modules parametrized by
  $W$}\\
\Points(W) &\text{equivalence classes of extended point modules parametrized by
  $W$}
\end{array}
\]

\section{Non-commutative symmetric algebras}
\label{ncsym}
\subsection{Generalities}
\label{ref-4.1-21}
We will consider the following particular case of a sheaf-$\ZZ$-algebra. Let
$\alpha:X\r S$, $\beta:Y\r S$ be smooth equidimensional maps of the
same relative dimension
and let $\Escr\in\shbimod_S(X-Y)$ be locally free. 

Define
\begin{equation}
\label{ref-4.1-22}
X_n=\begin{cases}
X&\text{if $n$ is even}\\
Y&\text{if $n$ is odd}
\end{cases}
\end{equation}
In a similar way we define
\begin{equation}
\label{ref-4.2-23}
\alpha_n=\begin{cases}
\alpha&\text{if $n$ is even}\\
\beta&\text{if $n$ is odd}
\end{cases}
\end{equation}
We define $\Escr^{\ast n}$ as in the introduction. I.e.
\begin{align*}
\Escr^{\ast n}=
\begin{cases}
\Escr^{\overbrace{\ast\cdots\ast}^n}&\text{if $n>0$}\\
\Escr&\text{if $n=0$}\\
{}^{\overbrace{\ast\cdots\ast}^{-n}}\Escr&\text{if $n<0$}
\end{cases}
\end{align*}
We then define
$\SS(\Escr)$ as the sheaf-$\ZZ$-sheaf-algebra generated by the
$\Escr^{\ast n}$ subject to the relations $i(\Oscr_{X_n})$. More precisely
\begin{align*}
\Ascr_{mn}&=
\begin{cases}
0&\text{if $n<m$}\\
\Oscr_{X_n} &\text{if $n=m$}\\
\Escr^{\ast m} & \text{if $n=m+1$}\\
\Escr^{\ast m}\otimes\cdots \otimes \Escr^{\ast n-1}/\\
\qquad (i(\Oscr_{X_{m}})\otimes \Escr^{\ast m+2}\otimes \cdots \otimes \Escr^{\ast n-1}
+\cdots\\\qquad\qquad\cdots+ \Escr^{\ast m}\otimes \cdots \otimes \Escr^{\ast n-3}\otimes i(\Oscr_{X_{n-2}}))
&\text{if $n\ge m+2$}
\end{cases}
\end{align*}

We say that
$\SS(\Escr)$ is a non-commutative symmetric algebra in \emph{standard
form}.

In the sequel it will sometimes be convenient to define more general
symmetric algebras. We will do so now and then we will show that these
more general symmetric algebras are equivalent to those in standard form.

Let $\alpha_n:X_n\r S$ be arbitrary smooth equidimensional maps of the
same relative dimension.
Assume that $(\Escr_n)_n$, $(\Qscr_n)_n$ are respectively a series of
locally free objects in $\shbimod(X_n-X_{n+1})$ and invertible objects
in $\shbimod(X_n-X_{n+2})$ which are non-degenerate subobjects of
$\Escr_n\otimes_{\Oscr_{X_{n+1}}}  \Escr_{n+1}$. 
 We then define
$\Ascr$ to be the $(X_n)_n$-sheaf-$\ZZ$-algebra generated by the
$\Escr_n$ subject to the relations $\Qscr_n$. Thus
$\Ascr_{nn}=\Oscr_{X_n}$, $\Ascr_{n,n+1}=\Escr_n$ and
$\Ascr_{n,n+2}=\Escr_n\otimes \Escr_{n+1}/\Qscr_n$, etc\dots. We will
call an algebra of the form $\Ascr$ a non-commutative symmetric
algebra.  We expect a
non-commutative symmetric algebra to have good homological properties
 but this has only been proved in the rank two case (see below).

Now let $X=X_0$, $\alpha=\alpha_0$, $Y=X_1$, $\beta=\alpha_1$ and
define $X'_n,\alpha'_n$ in the same way as $X_n,\alpha_n$ 
in \eqref{ref-4.1-22}\eqref{ref-4.2-23}. Thus
\[
X'_n=\begin{cases}
X_0=X&\text{if $n$ is even}\\
X_1=Y&\text{if $n$ is odd}
\end{cases}
\]
and
\[
\alpha'_n=\begin{cases}
\alpha_0=\alpha&\text{if $n$ is even}\\
\alpha_1=\beta&\text{if $n$ is odd}
\end{cases}
\]

Using \eqref{ref-3.5-13} we find:
\begin{align*}
\Escr_1&=\Escr_0^\ast \otimes_{\Oscr_X} \Qscr_0\\
\Escr_2&= \Qscr_0^{-1}\otimes_{\Oscr_X} \Escr^{\ast\ast}_0
 \otimes_{\Oscr_Y}\Qscr_1\\
\Escr_3&= \Qscr_1^{-1}\otimes_{\Oscr_Y} \Escr^{\ast\ast\ast}_0
\otimes_{\Oscr_X}\Qscr_0\otimes \Qscr_2
\end{align*}
Continuing we find that for $n\in\ZZ$ there exist invertible
$\Qscr_n'\in \shbimod(X'_n-X_n)$ such that 
\begin{equation}
\label{ref-4.3-24}
\Escr_n=\Qscr_n^{\prime -1}\otimes_{\Oscr_{X'_n}} \Escr_0^{\ast n} \otimes_{\Oscr_{X'_{n+1}}}
\Qscr'_{n+1}
\end{equation}
and 
\[
\Qscr_n=\Qscr_n^{\prime -1} \otimes_{\Oscr_{X'_n}}\Qscr'_{n+2}
\]
 The inclusion
\[
\Qscr_n\hookrightarrow \Escr_n\otimes_{\Oscr_{X_{n+1}}} \Escr_{n+1}
\]
becomes an inclusion
\[
\Qscr_n^{\prime -1} \otimes_{\Oscr_{X'_n}}\Qscr'_{n+2}\hookrightarrow
\Qscr_n^{\prime -1}\otimes_{\Oscr_{X'_n}} \Escr_0^{\ast n}
\otimes_{\Oscr_{X'_{n+1}}} \Escr_0^{\ast(n+1)} \otimes_{\Oscr_{X'_{n+2}}} \Qscr'_{n+2}
\]
and it is easy to see that this inclusion is derived from the canonical
inclusion 
$i_n:\Oscr_{X'_n}\r \Escr_0^{\ast n}\otimes_{\Oscr_{X'_{n+1}}} \Escr^{\ast(n+1)}$.  

Thus we have shown that  every non-commutative symmetric
algebra is obtained from one in standard form by twisting (see \S
\ref{ref-3.2-14}).

We will say that $\Ascr$ is a non-commutative symmetric algebra of
rank $r$ if $\Escr_0$ has rank $r$ on both sides. From lemma
\ref{ref-3.1.9-11}  together with \eqref{ref-4.3-24} we then obtain
that all $\Escr_n$ have rank $r$ on both sides.
\subsection{Relation with the  definition from \cite{VdB11,Pat1,Pat2} }
\label{ref-4.2-25}
Let $X$ be a scheme and let $\Escr\subset\shbimod_S(X-X)$ be locally
free. Let $\Qscr\in \Escr\otimes_{\Oscr_X}\Escr$ be a non-degenerate
invertible subobject
and let  $\Hscr=T_X(\Escr)/(\Qscr)$. The following lemma makes
the connection between $\Hscr$ and $\SS(\Escr)$.
\begin{lemmas}
We have $\Gr(\Hscr)\cong \Gr(\SS(\Escr))$. 
\end{lemmas}
\begin{proof}
If $\Ascr$ is a sheaf-$\ZZ$-graded algebra on $X$ then we
define the $\ZZ$-graded sheaf algebra $\check{\Ascr}$ by
\begin{equation}
\label{ref-4.4-26}
\check{\Ascr}_{ij}=\Ascr_{j-i}
\end{equation}
It is clear that we have $\Gr(\check{\Ascr})=\Gr(\Ascr)$. Furthermore it
is also clear that $\check{\Ascr}$ is a non-commutative symmetric
algebra with $\Escr_i=\Escr$ and $\Qscr_i=\Qscr$ for all $i$. Since
such a non-commutative symmetric algebra is obtained by twisting from
$\SS(\Escr)$ we are done.
\end{proof}
\subsection{Point modules over non-commutative symmetric algebras or
  rank two}
We let the notations be as in the previous sections but we assume in
addition that $\Ascr$ has rank two. We start with the
following result. 
\begin{propositions} 
\label{ref-4.3.1-27}
\label{ref-4.3.1-28} Assume that $P_{[m:m+1]}$ is a $[m:m+1]$-truncated point
  module over $\Ascr$. Then there exist unique (up to
  isomorphism) $[m-1:m+1]$ and $[m:m+2]$-truncated point modules $P_{[m-1:m+1]}$ and $P_{[m:m+2]}$
  whose restriction is equal to $P_{[m:m+1]}$. 
\end{propositions}
\begin{proof} Both claims are similar so we only consider the second
  one. Since we may shift $\Ascr$ we may without loss of generality assume that
  $m=0$. In that case $P$ is  described by a triple $(P_0,P_1,\phi)$ where
  $P_0\in\coh(X_0)$, $P_1\in\coh(X_1)$ are locally free of rank one
  over $S$ and $\phi:P_0\otimes_{\Oscr_X}
  \Escr_0\r P_1$ is a surjective map. We have to extend this triple to
  a quintuple $(P_0,P_1,P_2,\phi,\psi)$ where $P_2\in\coh(X_2)$ is
  also locally free of rank one over $S$ and
  $\psi:P_1\otimes_{\Oscr_{X_1}} \Escr_1\r P_2$ is another surjective map. The
  entries in such a quintuple are not arbitrary since the relation
  $\Qscr_0$ has to be satisfied. 
To clarify this restriction we note that point modules and truncated
point modules are preserved under twisting (see
\S\ref{ref-3.2-14}). Hence we 
may without loss of generality assume that $\Ascr$ is in standard
form, i.e. $\Ascr=\SS(\Escr)$ for some sheaf-bimodule $\Escr$ which is locally
free of rank two on both sides. 

In order for $(P_0,P_1,P_2,\phi,\psi)$ to define an object in
$\Gr(\Ascr_{[0:2]})$ 
module we need that the composition $P_0 \r P_0\otimes_{\Oscr_{X_0}}
\Escr\otimes_{\Oscr_{X_1}} \Escr^\ast \xrightarrow{\phi\otimes \Escr^\ast}
P_1 \otimes_{\Oscr_{X_1}} \Escr^\ast \xrightarrow{\psi} P_2$ is
equal to zero since this composition represents the action of
$\Qscr_0$. From lemma \ref{ref-4.3.2-30} below it follows that this composition
may be described in the following alternative way:
\begin{equation}
\label{ref-4.5-29}
P_0\xrightarrow{\phi^\ast} P_1\otimes_{\Oscr_{X_1}} \Escr^\ast
\xrightarrow{\psi} P_2
\end{equation}
where $\phi^\ast$ is obtained from $\phi$ by adjointness. Thus the
pair $(\psi,P_2)$ is a quotient of  $\coker \phi^\ast$. If we now show
that $\coker \phi^\ast$ is itself locally free of rank one then we are
done. This last fact follows from lemma \ref{ref-4.3.4-31} below.
\end{proof}
\begin{lemmas} 
\label{ref-4.3.2-30}
Assume that $(L,R)$ is a pair of adjoint functors and
  assume that we have objects $A,B$, together with a map $\phi:LA\r B$.
 Then the composition $A\r RL A\xrightarrow{R\phi} RB$ is equal 
to 
  $\phi^\ast:A\r RB$.
\end{lemmas}
\begin{proof}
This is standard.
\end{proof}
\begin{lemmas}
\label{rklemma2}
Let $\Escr\in \shbimod_S(X-Y)$ be
  locally free on 
the left
and let $\Fscr$ be a coherent $\Oscr_X$-module which  is locally free
over $S$. Then $\Fscr\otimes_{\Oscr_X} \Escr$ is also 
locally free over $S$. If $\Escr$ has constant rank $m$ on the left and similarly if
the $S$-rank of $\Fscr$ is constant and equal to $n$ then the $S$-rank of $\Fscr\otimes_{\Oscr_X}\Escr$
is constant as well and equal to $mn$. 
\end{lemmas}
\begin{proof} This is a direct consequence of Lemma \ref{rklemma} if we view $\Fscr$ as
an $S-X$-bimodule.
\end{proof}
\begin{lemmas}
\label{ref-4.3.4-31}
Let $\alpha:X\r S$, $\beta:Y\r S$ be smooth equidimensional maps of
the same relative dimension. 
 Let $\Escr\in \shbimod_S(X-Y)$ be
  locally free of rank two on both sides.  Assume that we have objects
  $P_0\in \coh(X)$, $P_1\in \coh(Y)$ which are locally free of rank one over $S$,
  together with a surjective map $\phi:P_0\otimes_{\Oscr_X} \Escr\r
  P_1$. Then the adjoint map $\phi^\ast:P_0\r 
  P_1\otimes_{\Oscr_Y} \Escr^\ast$ is injective and has a cokernel which is
  locally free of rank one over $S$.
\end{lemmas}
\begin{proof}
Using lemma \ref{ref-3.1.5-4} it suffices to prove this in the case
that $S=\Spec k$. But then it is sufficient to show that $\phi^\ast$
is not zero (as $P_1\otimes_{\Oscr_Y} \Escr^\ast$ has rank two by Lemma \ref{rklemma2}). Since $\phi$ is not zero this is clear.
\end{proof}
Using the bijections exhibited in Proposition \ref{ref-4.3.1-28} together with
the fact that the relations in $\Ascr$ have degree two we now easily
obtain:
\begin{theorems} 
\label{ref-4.3.5-32}
The sets of extended point modules, $m$-shifted point
modules and $[m:n]$-truncated point modules for $n\ge m+1$ over $\Ascr$ are all in
bijection.  These bijections are given by the appropriate restriction functors.
\end{theorems}
\begin{corollarys}
The functors $\Points_\Ascr$, $\Points_{m,\Ascr}$ and
$\Points_{m,n,\Ascr}$ (for $n\ge m+1$) are all
naturally equivalent.
\end{corollarys}

It follows from the proof of Proposition \eqref{ref-4.3.1-27}
(see \eqref{ref-4.5-29})
that if $P$ is an extended point module over $\Ascr$ then there are exact sequences
 on $X_{j+2}$ 
\begin{equation}
\label{ref-4.6-33}
0\r P_j\otimes_{\Oscr_{X_{j}}}  \Qscr_j \r P_{j+1}
\otimes_{\Oscr_{X_{j+1}}}
\Escr_{j+1}  \r P_{j+2}\r 0
\end{equation}
In fact this was only shown if $\Ascr$ is in standard form, but the
general case follows by twisting. Now write $P_j$ in the usual form
$\zeta_{j\ast}(Q_j)$ where $\zeta_j$ is a section of $\alpha_j$ and $Q_j\in
\Pic(S)$. Then applying $\alpha_{j+2\ast}$  to \eqref{ref-4.6-33} we obtain
an exact sequence on $S$ 
\[
0\r Q_j\otimes_{\Oscr_S} \zeta^\ast_j \pr_{1\ast}(\Qscr_j)\r
Q_{j+1}\otimes_{\Oscr_S} \zeta^\ast_{j+1} \pr_{1\ast}(\Escr_{j+1})\r
Q_{j+2}
\r
0
\]
Put $\Bscr=\Bscr(P)$. Tensoring the previous exact sequence on the left with $Q^{-1}_i$ yields an exact sequence
\begin{equation}
\label{ref-4.7-34}
0\r \Bscr_{ij}\otimes_{\Oscr_S} \zeta^\ast_j \pr_{1\ast}(\Qscr_j)\r
\Bscr_{ij+1}\otimes_{\Oscr_S} \zeta^\ast_{j+1} \pr_{1\ast}(\Escr_{j+1})\r
\Bscr_{ij+2}
\r
0
\end{equation}

By dualizing \eqref{ref-4.6-33}, tensoring on the left with $\Qscr_j$, applying a
suitable variant of \eqref{ref-3.5-13}, applying $\alpha_{j\ast}$,
tensoring with $Q_k$ and finally changing indices we obtain the
following analogous exact sequence 
\begin{equation}
\label{ref-4.8-35}
0\r \zeta_{i+2}^\ast\pr_{2\ast}(\Qscr_i)\otimes_{\Oscr_S} \Bscr_{i+2j}
\r
\zeta_{i+1}^\ast\pr_{2\ast}(\Escr_{i})\otimes_{\Oscr_S}\Bscr_{i+1j}
\r
\Bscr_{ij}
\r
0
\end{equation}

\subsection{Projective bundles associated to quasi-coherent sheaves}
 If $Z$ is a scheme and $\Uscr$ is a coherent
sheaf on $Z$ then we define $\PP_Z(\Uscr)=\underline{\Proj} \,S_Z\Uscr$
where $S_Z\Uscr=\oplus_n S^n_Z\Uscr$ denotes the symmetric algebra of
$\Uscr$.  On $E=\PP_Z(\Uscr)$ there is a canonical line bundle denoted
by $\Oscr(1)$ or $\Oscr_E(1)$
which corresponds to $(S_Z\Uscr)(1)$.

If $W$ is an arbitrary scheme and $\chi$ is a $W$-point of
$\PP_Z(\Uscr)$ then $\chi$ defines a pair $(\chi',\Lscr)$ where
$\chi'$ is the composition $W\xrightarrow{\chi}\PP_Z(\Uscr)\r Z$ and
$\Lscr\in\Pic(W)$ is given by $\chi^\ast(\Oscr(1))$. Clearly $\Lscr$
is a quotient of $\chi^{\prime \ast}(\Uscr)$. It is standard
that conversely every pair $(\chi',\Lscr)$ where $\chi'$ is a map $W\r
Z$ and $\Lscr\in\Pic(W)$ is a quotient of $\chi^{\prime\ast}(\Uscr)$
corresponds to a unique $\chi:W\r \PP_Z(\Uscr)$.

We will use the following result in the following sections.
\begin{lemmas} 
\label{ref-4.4.1-36}
Let $x\in Z$ and let $m_x\subset \Oscr_{Z,x}$ be the
  maximal ideal. Then the scheme-theoretic closed fiber of  $x$ in
  $\PP_Z(\Uscr)$ is equal to $\PP_{k(x)} (\Uscr_x/m_x \Uscr_x)$. In
  particular it is equal to some $\PP^n_{k(x)}$.
\end{lemmas}
Here is a somewhat more specialized result.
\begin{propositions}
\label{ref-4.4.2-37}
Assume that $\beta:Z\r X$ is a  map of schemes and assume that
$\Escr\in \coh(Z)$ is coherent over $X$. Then the
obvious map $o:\PP_Z(\Escr)\r \PP_X(\beta_\ast\Escr)$ is a closed
immersion.  If $X$ is a smooth 
connected curve over $k$  and $\Escr$ is locally free of rank two over
$X$ then $o$ is either surjective
 or else its image is a divisor.
\end{propositions}
\begin{proof}
All claims are local on $X$ so we may and we will assume that $X=\Spec R$ is
affine. 
In addition we may replace $Z$ by
the scheme-theoretic support of $\Escr$, i.e. we may assume that
$\beta$ is finite. It follows that $Z$ is also
affine, say $Z=\Spec T$. Therefore $\Escr$ is obtained from a finitely generated $T$
module $E$ and $\PP_Z(\Escr)=\Proj S_T(E)$, $\PP_X(\Escr)=\Proj
S_R(E)$. 
The map $o$ is obtained from the obvious map $S_R(E)\r
S_T(E)$. 

To prove that $o$ is a closed immersion we simply remark that $S_T(E)\r
S_R(E)$ is surjective in degree $\ge 1$.

Now we make the additional hypotheses on our data, i.e.   $X$ is a smooth 
connected curve over $k$  and $\Escr$ is locally free of rank two over
$X$. To prove our claim we may now make the additional simplifying
assumption that $X=\Spec R$ where $R$ is a
discrete valuation ring.

The fact that $E$ is Cohen-Macaulay
implies that $T$ has no embedded components. So $T$ is free of rank
one or two over $R$ and  $R$ embeds in $T$.

If $T$ is free of rank one then $T=R$ and hence $o$ is an
isomorphism. So assume that $T$ has rank two. 
Thus
$T=R[z]$ where $z$ satisfies 
a monic quadratic equation over $R$.

We now have to show that the kernel $K$ of  $S_R(E)\r
S_T(E)$ is generated by one element. Let $E=Rx+Ry$. Then $K$ is
generated by $(z\cdot x)x-x(z\cdot x)$, $(z\cdot y)x-y(z\cdot x)$ and
$(z\cdot y)y-y(z\cdot y)$.
Write $z\cdot x=ax+by$, $z\cdot y=cx+dy$ with $a,b,c,d\in
R$. Then
\begin{align*}
(z\cdot x)x-x(z\cdot x)&=byx-bxy=0\\
(z\cdot y)x-y(z\cdot x)&=cxx+dyx-ayx-byy=cx^2+(d-a)xy-by^2\\
(z\cdot y)y-y(z\cdot y)&=dxy-dyx=0
\end{align*}
Thus $K$ is indeed generated by a single quadratic element.
\end{proof}
\begin{remarks} 
\label{ref-4.4.3-38}
The preceding result is false if $X$ is not a
  curve.

Consider the following example~: $Z=X\times Y$ with $X=Y=A^2$, $\Delta\subset
X\times Y$ is the diagonal and $\Gamma$ is the graph of $(x,y)\mapsto
(-x,-y)$. Let $\Escr=\Oscr_\Delta\oplus \Oscr_\Gamma$. Counting
dimensions of fibers we see that $\PP_{X\times Y}(\Escr)$ has
dimension $2$.

Clearly
$\PP_{X\times Y}(\Escr)$ contains two closed subsets respectively given by
$\PP_{X\times Y}(\Oscr_{\Delta})=\Delta$ and
$\PP_{X\times Y}(\Oscr_{\Gamma})=\Gamma$ which must be irreducible
components since they also have dimension $2$. Furthermore outside the
point $(o,o)\in X\times Y$ the map $\Delta\coprod\Gamma \r \PP_{X\times
  Y}(\Escr)$ is an isomorphism.  However  the fiber $F$ of
$(o,o)$ in $\PP_{X\times
  Y}(\Escr)$ is $\PP^1$ whereas $\Delta\coprod\Gamma$ gives us at most
two points. 

Thus $F$ must be contained in an additional irreducible component. If
this irreducible component is not $F$ itself then it must contain some
points of $\PP_{X\times
  Y}(\Escr)$ not above $(o,o)$. But then $F$ must be equal to $\Delta$
or $\Gamma$ which is a contradiction. It follows that $\PP_{X\times
  Y}(\Escr)$ is not equidimensional and in particular it cannot be a
divisor in $\PP_X(\pr_1\Escr)$.

 The problem with this example is that the support $\Delta\cup\Gamma$ of $\Escr$  is not Cohen-Macaulay.
\end{remarks}
\subsection{Representability of the point functor}
The following result has been proved by Adam Nyman \cite{Nyman1}. We
reproduce the proof since we need the exact nature of the
isomorphisms involved.
\begin{theorems} The functor $\Points_\Ascr$ is representable by
  $\PP_{X\times_S Y}(\Escr_0)$.
\end{theorems}
\begin{proof}
In view of the above discussion it is clearly sufficient to prove this
for $\Points_{0,1,\Ascr}$. We will start by giving an alternative
description of $\Points_{0,1,\Ascr}(S)$. Without loss of generality we may  assume that $\Ascr=\SS(\Escr)$.

An object in $\Points_{0,1,\Ascr}(S)$ has a unique representative of the
form $(P_0,P_1,\phi)$ where $\alpha_{0,\ast}(\Pscr_0)=\Oscr_S$ and
$\phi:P_0\otimes_{\Oscr_X} \Escr\r P_1$ is an
epimorphism. There exist sections $\zeta_0,\zeta_1$ of $\alpha,\beta$
and an element $Q_1$ of $\Pic(S)$ 
such that $P_0=\zeta_{0,\ast} (\Oscr_S)$ and $P_1=\zeta_{1,\ast}(Q_1)$.

 According to lemma \ref{ref-3.6-19} $\phi$ corresponds to an
epimorphism $\phi':\Escr\r P_0\boxtimes_S P_1$ and furthermore
$P_0\boxtimes_S P_1=(\zeta_0,\zeta_1)_\ast(Q_0)$.  Since
$(\zeta_0,\zeta_1)_\ast(Q_0)$ contains all information to reconstruct
$\zeta_0,\zeta_1$ and $Q_0$ we conclude that $\Points_{\Ascr,0,1}(S)$
is in one-one correspondence with the set of quotients of $\Escr$ on $X\times_S Y$
which are of rank one over $S$.

If we apply this the discussion before the statement of the theorem
with  $Z=X\times_S Y$, $\Uscr=\Escr$, $W=S$ then we
find 
\[
\Points_{\Ascr,0,1}(S)=\Hom_{\Sch}(S,\PP_{X\times_S Y}(\Escr))
\]
Since this bijection is obviously compatible with base extension we
find that the functor   $\Points_{\Ascr,0,1}$ is represented by
$\PP_{X\times_S Y}(\Escr)$. This finishes the proof.
\end{proof}

\section{Properties of the universal point algebra}
From now on we assume that our base scheme $S$ is $\Spec k$ and
therefore we will omit $S$ from the notations. Otherwise the notations will be as in the previous section.
\subsection{A vanishing result} 
\begin{theorems}
\label{ref-5.1.1-39}
Let $s:E\r \bar{E}$ be a projective map  of relative dimension one.
Let $\Lscr\in \coh(E)$ and assume that the restriction to every fiber
of $\Lscr$ 
is generated by global sections and has vanishing higher cohomology. Then $R^i
s_\ast\Lscr=0$ for $i>0$ and the canonical map $s^\ast s_\ast
\Lscr\r \Lscr$ is a surjective.
\end{theorems}
\begin{proof}
This is not an immediate consequence of semi-continuity since we are not
assuming that $\Lscr$ is flat over $\bar{E}$. 

We use the theorem on formal functions. For $y\in \bar{E}$ let
$E_n=E\times_{\bar{E}} \Spec \Oscr_{\bar{E},y}/m_y^n$ where $m_y$ is
the maximal ideal corresponding to $y$. In addition let $\Lscr_n$ be
the restriction of $\Lscr$ to $E_n$.
Then one has \cite[Thm III.11.1]{H}
\[
(R^is_\ast\Lscr)\,\hat{}_y= \projlim_n H^i(E_n,\Lscr_n)
\]
Thus in order to show that   $R^i
s_\ast(\Lscr)=0$ for $i>0$ it is sufficient to show that
\begin{itemize}
\item[($H1_n$)]
$H^i(E_n,\Lscr_n)=0$ 
\end{itemize} 
for all $y$ and all $n$.

Similarly it is easy to see that for   $s^\ast s_\ast
\Lscr\r \Lscr$ to be surjective  it is sufficient that the condition
\begin{itemize}
\item[($H2_n$)]
$\Gamma(E_n,\Lscr_n)\otimes_k \Oscr_{E_n}\r \Lscr_n$ is surjective
\end{itemize}
holds for all $y$ and all $n$. 

Our proof will be by induction on $n$.  It
follows from the hypotheses that $(H1_1)$ and $(H2_1)$ are
satisfied. 

Assume now that $(H1_n)$ and $(H2_n)$ are satisfied.  We have an exact
sequence
\[ 
m^n_y /m^{n+1}_y\otimes_k \Lscr_1 \r  \Lscr_{n+1} \r \Lscr_n
\r 0
\]
Thus $\Fscr=\ker (\Lscr_{n+1}\r\Lscr_n)$ is  the quotient of a sheaf
with vanishing higher cohomology, and since we are in dimension one it
follows that $\Fscr$ itself has vanishing higher cohomology.

 Thus it follows
that 
\[
0\r H^0(E_1,\Fscr)\r H^0(E_{n+1},\Lscr_{n+1})\r H^0(E_n,\Lscr_n)\r 0
\]
is exact, and furthermore the induction hypotheses imply that
$H^i(E_{n+1},\Lscr_{n+1})=0$ for $i>0$. So this proves $(H2_{n+1})$. 

In order to prove $(H1_{n+1})$ we use the following 
commutative diagram with exact rows:
\[
\begin{CD}
0@>>>\Fscr @>>>  \Lscr_{n+1} @>>> \Lscr_n
@>>> 0 \\
@. @AAA @AAA @AAA @.\\
0@>>> H^0(E_1,\Fscr)\otimes_k \Oscr_E@>>>  H^0(X,\Lscr_{n+1})\otimes_k \Oscr_E
@>>>H^0(X,\Lscr_{n})\otimes_k \Oscr_E  @>>> 0
\end{CD}
\]
Since the outermost vertical maps are surjective the same holds for
the middle one. This proves $(H1_{n+1})$.
\end{proof}
\subsection{The case of non-commutative symmetric algebras}
\label{ref-5.2-40}
In this section the notations are as before. In particular
$\Ascr$ is a non-commutative symmetric algebra of rank two  over
$\Xi=(X_i)_{i\in\ZZ}$ (see \S\ref{ncsym}).  As usual we put $\Escr_i=\Ascr_{i,i+1}$.
By definition $\Escr_i$ has rank two on both sides. 

Put $E^j=\PP_{X_j\times X_{j+1}} (\Escr_j)$. 
 Since $E^j$ represents
$\Points_{\Ascr}$, there is a universal extended point $P^j$ over
$\Ascr_{E^j}$. We now let $\Bscr^j=\Bscr(P^j)$,
 be the associated
sheaf-$\ZZ$-algebras and we aim to study these in more detail. As
above let $\zeta^j_i:E^j \r X_{i,E^j}$ be the support of $P^j_i$. We may
write $\zeta^j_i$ as a pair $(\mu^j_i,\Id_{E^j})$ where $\mu^j_i$ is a
map from $E^j$ to $X_i$. Again as above we write
$P^j_i=\zeta^j_{i,\ast}(Q^j_i)$ for $Q^j_i\in \Pic(E^j)$.  We will
also write $\alpha^j_i:X_{i,E^j}\r E^j$ for the map obtained by base
extension from $\alpha_i:X_i\r \Spec k$.

Our first observation is that since the $E^j$ all represent the same
functor there must exist isomorphisms $\theta^j:E^{j+1}\r E^j$ and
objects $L^j\in\Pic(E^j)$ such that 
\[
P^{j+1}_i=\alpha_i^{j+1\ast}(L^{j+1})\otimes_{\Oscr_{X_{i,E^{j+1}}}}
\theta^{j\ast}(P^j_i) 
\]
This may be rewritten as $\mu^{j+1}_i=\mu^j_i \theta^j$ and
$Q^{j+1}_i=L^{j+1} \otimes_{\Oscr_{E^{j+1}}} \theta^{j\ast} Q^j_i$ from which
we deduce 
\[
 {\Bscr}^{j+1}_{mn}= \theta^{j\ast}({\Bscr}^{j}_{mn})
\]
In the sequel we will define $\theta^{jl}:E^j \r E^l$ as the
composition $\theta^l \theta^{l+1} \cdots \theta^{j-1}$ if $j\ge l$
and by a similar formula if $j<l$. Thus we find 
\[
\mu^j_i   =  \mu^l_i \theta^{jl}
\]
\and
\[
{\Bscr}^{j}_{mn}=\theta^{jl\ast} {\Bscr}^l_{mn}
\]
From the proof that $E^m$ represents $\Points_\Ascr$ it follows that
${\Bscr}^m_{m,m+1}=\Oscr_{E^m}(1)$ 
 and $(\mu^m_{m}, \mu^m_{m+1})$ is
the projection map $E^m=\PP_{X_m\times X_{m+1}}(\Escr_m)\r X_m\times X_{m+1}$. This allows us to
describe ${\Bscr}_{mn}^i$ in terms of the $\Oscr_{E^j}(1)$ and the
isomorphisms $\theta^{pq}$.

Let $E^j\xrightarrow{s^j_i}\bar{E}^j_i\xrightarrow{\bar{\mu}^j_i}X_i$ be the
Stein factorization of $\mu^j_i$. To understand these factorizations
let us first consider $\mu^j_j$  and $\mu^j_{j+1}$ which together represent the
canonical map $E^j
\r
X_j\times X_{j+1}$. As an intermediate step consider the Stein
factorization $E^j
\r G^j\r
X_j\times X_{j+1}$ of this last map.   By construction \cite[Cor. III.11.5]{H} $G^j$ is finite over the scheme
theoretic image $Z^j$ of $E^j$ in $X_j\times X_{j+1}$. Since $Z^j$
is finite over both $X_j$ and $X_{j+1}$ we obtain from the
construction of $\bar{E}^j$ \cite[Cor. III.11.5]{H} that $E^j\r
 G^j\r X_{j}$ and $E^j\r G^j
\r
 X_{j+1}$ are the Stein factorizations of
respectively $\mu^j_j:E^j\r X_j$ and $\mu^j_{j+1}:E^j\r X_{j+1}$.  In
particular we obtain $\bar{E}^j_j=\bar{E}^j_{j+1}$ and
$s^j_j=s^{j}_{j+1}$.   Now using the fact that Stein factorizations are
(obviously) 
compatible with isomorphisms we obtain from this by applying suitable
$\theta^{pq}$ that $\bar{E}^p_j=\bar{E}^p_{j+1}$ and
$s^p_j=s^{p}_{j+1}$ for all $p$. Thus $\bar{E}^p_j$ and $s^p_j$ are
independent of $j$ and we may write $\bar{E}^p_j=\bar{E}^p$, $s^p_j=s^p$.
Thus the result of this discussion is that we have commutative
diagrams:
\begin{equation}
\label{ref-5.1-41}
\begin{CD}
E^p @>\theta^{pq}>> E^q\\
@V s^p VV   @VV s^q V\\
\bar{E}^p @>\bar{\theta}^{pq} >> \bar{E}^q\\
@V \bar{\mu}^p_i VV  @VV \bar{\mu}^q_i V\\
X_i @= X_i
\end{CD}
\end{equation}

Now we investigate the scheme-theoretic closed fibers of $s^j$. 

By lemma \ref{ref-4.4.1-36} the scheme-theoretic fibers of 
$E^j\r Z^j$ are either points or $\PP^1$'s and hence in particular
they are connected.
The fibers of  $E^j\r \bar{E}^j$  are also connected by the properties
of the Stein factorization.
Hence it
follows that the map $\bar{E}^j\r Z^j$ is settheoretically a
bijection. In particular $s^j$ and $E^j \r Z^j$ have the same closed
fibers. We conclude that the fibers if $s^j$ are either points or $\PP^1$'s.

Now let $\Omega$ be the constant system of schemes $(E^0)_{i\in\ZZ}$
and let $\mu^0=(\mu^0_i)_i$. From Corollary \ref{ref-3.4.4-20}
it follows that $\mu^0$ satisfies
condition (C). We can now prove the following technical result which will be used
below in the proof that a non-commutative symmetric algebra is noetherian (see \S\ref{secnoetherian} below).
\begin{theorems}
\label{ref-5.2.1-42}
${\Bscr}^0_{\ge 0}$ is ample for $\mu^0$.
\end{theorems}
\begin{proof}
 Since ${\Bscr}^0$  is strongly
graded and ${\Bscr}^0_{00}=\Oscr_{E^0}$ it is clear that
${\Bscr}^0$ and hence ${\Bscr}^0_{\ge 0}$  is noetherian. So we need
only  verify the conditions
2. and 3. from Theorem \ref{ref-3.3.2-17}. Let $\Mscr\in\coh(E^0)$.

We compute 
\begin{align*}
\Mscr\otimes_{\Oscr_{E^0}}{\Bscr}^0_{ij}&=
\Mscr\otimes_{\Oscr_{E^0}}{\Bscr}^0_{i,i+1}
\otimes_{\Oscr_{E^0}}{\Bscr}^0_{i+1,i+2}\otimes_{\Oscr_{E^0}}
\cdots \otimes_{\Oscr_{E^0}} {\Bscr}^0_{j-1,j}\\
&=\Mscr\otimes_{\Oscr_{E^0}} \theta^{i0}_\ast
(\Bscr^i_{i,i+1})
\otimes_{\Oscr_{E^0}}\cdots
\otimes_{\Oscr_{E^0}}\theta^{j-1,0}_\ast(\Bscr^{j-1}_{j-1,j})\\
&=\Mscr\otimes_{\Oscr_{E^0}} \theta^{i0}_\ast
(\Oscr_{E^i}(1))
\otimes_{\Oscr_{E^0}}\cdots
\otimes_{\Oscr_{E^0}}\theta^{j-1,0}_\ast(\Oscr_{E^{j-1}}(1))
\end{align*}
Since
\[
R\mu^0_{j\ast} (\Mscr\otimes_{\Oscr_{E^0}}{\Bscr}^0_{ij})=
\bar{\mu}^0_{j\ast} R
s^0_\ast(\Mscr\otimes_{\Oscr_{E^0}}{\Bscr}^0_{ij})
\]
and $\bar{\mu}^0_{j}$ is finite, it is sufficient to prove the analogues of 2. and 3. in Theorem
\ref{ref-3.3.2-17} for $R^i
s^0_\ast(\Mscr\otimes_{\Oscr_{E^0}}{\Bscr}^0_{ij})$. According to
Theorem \ref{ref-5.1.1-39} we have to show that 
$\Mscr\otimes_{\Oscr_{E^0}}{\Bscr}^0_{ij}$ when restricted to the
$\PP^1$ fibers of $s^0$ becomes eventually generated by global sections. This follows from
the fact that according to \eqref{ref-5.1-41} the
$\PP^1$-fibers are preserved under the $\theta$'s and the fact that
 $\Oscr_{E^m}(1)$ when restricted to a $\PP^1$-fiber of $s^m$ is equal
 to $\Oscr_{\PP^1}(1)$. 
\end{proof}
\section{On the structure of non-commutative symmetric algebras of
  rank two}
In this section the notations are the same as in the previous ones. 
\subsection{Ranks and exact sequences}
Let $e_i\in\Gamma(X_n,\Ascr_{nn})=\Gamma(X_n,\Oscr_{X_n})$ be the section corresponding to $1$.
The structure of the relations in $\Ascr$ implies that there is an exact
sequence of $\Oscr_{X_m}-\Ascr$ sheaf-bimodules given by
\begin{equation}
\label{ref-6.1-43}
\Qscr_m\otimes_{\Oscr_{X_{m+2}}}e_{m+2}\Ascr\r \Escr_{m}\otimes_{\Oscr_{X_{m+1}}}
e_{m+1}\Ascr \r e_m\Ascr\r 0
\end{equation}
We will show below that this exact sequence is exact on the left. 

The following proposition is proved in the same way as Proposition
\ref{ref-4.3.1-28} and Theorem \ref{ref-4.3.5-32}.
\begin{propositions}
\label{ref-6.1.1-44}
Assume that $Q_{[0:n]}$ is an object in $\Gr(\Ascr_{[0:n]})$  with
the following properties
\begin{enumerate}
\item
$(Q_{[0:n]})_i\neq 0$ for all $i\in \{0,\ldots,n\}$.
\item
$Q_{[0:n]}$ is generated in degree zero.
\item  $(Q_{[0:n]})_0$ and $(Q_{[0:n]})_1$ have finite length and
$\dim H^0(Q_{[0:n]})_0=\dim H^0(Q_{[0:n]})_1=1$.
\end{enumerate}
Then $Q_{[0:n]}$ is a $[0:n]$-truncated point module. Similarly if $Q$ is an
object in $\Gr(\Ascr)$ satisfying  suitable analogues  of (1-3) then $Q$ is
a point module. 
\end{propositions}
From the fact that a point module
is uniquely determined by its restriction to $\Ascr_{[0:1]}$ one
obtains that if $k$ is algebraically closed then for every $x\in X$ there is at least one point module $P$
such that $P_0=\Oscr_x$. 

Now we will consider line-modules. For $x$ a rational point in $X_m$ we define
$L_{m,x}=\Oscr_x\otimes_{\Oscr_{X_m}}e_m\Ascr$.  For simplicity we
write $L_x$ for $L_{0,x}$.

If $P$ is a point module then we have
\[
\Hom_{\Ascr}(L_x,P)=\Hom_{\Oscr_X}(\Oscr_x,P)
=\begin{cases}
k & \text{if $P_0=\Oscr_x$}\\
0 &\text{otherwise}
\end{cases}
\]
Thus it follows that if $k$ is algebraically closed then  every $L_x$ maps onto at least one
point-module. In the same way one sees that $L_{m,x}$ maps to an
$m$-shifted point module.

Let $L_x\r P$ be a surjective map to a point module and let $K$ be its
kernel. Since  $\length
(L_x)_1=2$ and $\length P_1=1$ we deduce that $K_1\cong \Oscr_y$ for
some $y\in X$. Thus there is a non-zero map $L_{1,y}\r K_1$. Since
$\coker (L_{1,y}\r L_x)$ has the same truncation to $\Ascr_{[0:1]}$ as
$P$ it
follows from Proposition \ref{ref-6.1.1-44} that  we have an exact sequence 
\begin{equation}
\label{ref-6.2-45}
L_{1,y}\r L_x\r P\r 0
\end{equation}
We will call this a standard exact sequence. A similar standard exact
sequence exists for $L_{m,x}$:
\begin{equation}
\label{ref-6.3-46}
L_{m+1,y}\r L_{m,x}\r P\r 0
\end{equation}
where $P$ is now an $m$-shifted point module. 

We can now prove the following result.
\begin{theorems} \label{new--1} We have
\begin{enumerate}
\item $\Ascr_{m,n}$ is locally free of rank $n-m+1$ on both sides.
\item The exact sequences \eqref{ref-6.1-43} and
  \eqref{ref-6.3-46} are exact on the left. 
\end{enumerate}
\end{theorems}
\begin{proof} Without loss of generality we may assume that $k$ is
  algebraically closed. As far as (1) is concerned we will only consider the
  left structure of $\Ascr$. The statement about the right structure
  follows by symmetry.

Assume that we have shown that $\Ascr_{m,n}$ is
  locally free on the left of rank $n-m+1$ for $n-m\le
  t$. We tensor \eqref{ref-6.1-43} on the left with
  $\Oscr_{x}$. Since $\length(\Oscr_x\otimes_{\Oscr_{X_m}}\Qscr_m)=1$
  and $\length(\Oscr_x\otimes_{\Oscr_{X_m}} \Escr_m)=2$ we obtain that
  $\Oscr_x\otimes_{\Oscr_{X_m} }\Qscr_m=\Oscr_{x'}$ and
  $\Oscr_x\otimes_{\Oscr_{X_m}} \Escr_m$ is an extension of
  $\Oscr_{x''}$ and $\Oscr_{x'''}$ for some $x',x'',x'''\in X$.

This yields
\begin{multline*}
\length (L_{m,x})_{m+t+1} \ge\\ \length(L_{m+1,x''})_{m+t+1}+
 \length (L_{m+1,x'''})_{m+t+1} - \length (L_{m+2,x'})_{m+t+1}=t+2
\end{multline*}
On the other hand we have from \eqref{ref-6.3-46} 
\begin{align*}
\length (L_{m,x})_{m+t+1}&\le 1+  \length (L_{m+1,y})_{m+t+1}\\
&=t+2
\end{align*}
Combining these two inequalities yields
$\length(L_{m,x})_{m+t+1}=t+2$ for all $m,x$. Since
$(L_{m,x})_{m+t+1}=\Oscr_{x}\otimes_{\Oscr_X} \Ascr_{m,m+t+1}$ this
yields that $\Ascr_{m,m+t+1}$ is locally free of rank $t+2$ on the left.

By induction we obtain the corresponding statement  for all $m,n$. From this we easily obtain that \eqref{ref-6.1-43} and
  \eqref{ref-6.3-46} are exact on the left. 
\end{proof}
\subsection{Two different types}
We need the following notation. Let $X=X'\bigcup X''$ and
$Y=Y'\bigcup Y''$ be disjoint unions of schemes  and let $p',p'':X',X''\r
X$, $q',q'':Y',Y''\r Y$ be the inclusion maps.
Assume $\Mscr'\in\ShBimod(X'-Y')$,  $\Mscr''\in\ShBimod(X''-Y'')$. Then we define $\Mscr'\boxplus
\Mscr''$ as $(p',q')_\ast(\Mscr')\oplus (p'',q'')_\ast(\Mscr'')$. 
We use a similar construction for sheaf-$\ZZ$-algebras. We leave the
obvious definitions to the reader.

We will now analyze the
$\Fscr\in\shbimod(X-Y)$ which are locally free of rank two on both
sides. As usual we assume that $X$, $Y$ are smooth of the same dimension and
equidimensional. 

Let $Z$ be the scheme theoretic support of $\Fscr$. Since $\Fscr$ is
Cohen-Macaulay, all components of $Z$ have the same dimension and
there are no embedded components.

  Assume that $Z$
has an irreducible component $Z'$ on which the restriction of $\Fscr$
has rank two (generically).  $Z'$ lies over connected components
$X'$ and $Y'$ of $X$ and $Y$. Let $X''$ and $Y''$ be the union of
the other connected components of $X$ and $Y$. Counting ranks we see
that there can be no other
irreducible components of $Z$ lying above $X'$ and $Y'$ and hence
$\Fscr=\Fscr'\boxplus \Fscr''$ where $\Fscr'\in \shbimod(X'-Y')$ and
$\Fscr''\in \shbimod(X''-Y'')$.

Let us return to $\Fscr'$. Since $Z'$ is integral and has degree one over $X'$ and $Y'$ and since
$X'$ and $Y'$ are furthermore integrally closed we obtain that $Z$ is
the graph of an isomorphism $\sigma:X\r Y$ and $\Fscr$ is a
vector bundle of rank two on $Z'$. 

 It is clear that $\SS(\Fscr)=\SS(\Fscr')\boxplus \SS(\Fscr'')$. A
similar decomposition then holds for every non-commutative symmetric
algebra by twisting. Furthermore we leave it to the
reader to check that $\Gr(\SS(\Fscr'))$ is equivalent to
$\Gr(S_{Z'}(\Fscr'))$ and hence corresponds to a commutative $\PP^1$-bundle.

To formalize this let us make the following definition.
\begin{definitions} Let $\Ascr$ be a non-commutative symmetric algebra
  of rank two and let $\Escr=\Ascr_{01}$. We say that $\Ascr$ is of
  Type I if $\Escr$ is a rank two bundle over the graph of an
  automorphism and we say that $\Ascr$ is of Type II if the restrictions
  of $\Escr$ to the irreducible components of its support all have
  rank one generically. 
\end{definitions}
Thus we have obtained the
following result.
\begin{propositions} Let $\Ascr$ be a non-commutative symmetric
  algebra of rank two. Then $\Ascr=\Ascr'\boxplus \Ascr''$ where
  $\Ascr'$ is of Type I and $\Ascr''$ is of type II. $\Gr(\Ascr')$ is
  equivalent to the category of graded modules over the symmetric
  algebra of a rank two vector bundle over a smooth scheme.
\end{propositions}
\begin{examples} The most basic example of type II symmetric algebra
  is obtained by embedding a smooth elliptic curve $C$ as a divisor of
  degree $(2,2)$ in $\PP^1\times \PP^1$ and letting
  $\Escr={}_u\Lscr_v$ where $\Lscr$ is a line bundle on $C$ and $(u,v):C\r \PP^1\times \PP^1$ denotes the embedding.  Such non-commutative symmetric algebras appeared naturally
in~\cite{VdB11} and provided one of the motivations for writing the current paper. 
\end{examples}
\subsection{Non-commutative symmetric algebras of rank two are noetherian}
\label{secnoetherian}
Since to prove $\Ascr$ is noetherian we may treat the
cases of Type I and Type II individually, and since the Type I case is
easy we assume throughout that $\Ascr$ is of Type II.

From Theorems \ref{ref-5.2.1-42} and \ref{ref-3.3.2-17} we obtain that
$\mu_\ast{\Bscr^0}_{\ge 0}$ is noetherian.  Furthermore by
construction there is a map $\Ascr\r \mu^0_\ast{\Bscr^0}_{\ge 0}$. We
would like to use this map in order to analyze $\Ascr$. However the
analysis is complicated by the fact that $E^0$ may have components of
different dimensions if $\dim X_n>1$ (see Remark \ref{ref-4.4.3-38}).

Therefore we will use the following trick. We will let $F^j$ be the
union of all components in $E^j$ which are of maximal dimension and we let
$t^j:F^j\r E^j$ be the inclusion map.  It is clear that 
with  $\theta^{jl}$ restricts to a map $F^j\r F^l$ which we will also
denote by $\theta^{jl}$.

Let
$\Cscr_{mn}=t^\ast({\Bscr^0}_{mn})$. Then
$\Cscr=\oplus_{m\le n}\Cscr_{mn}$ is a $\ZZ$-algebra on $F^0$. 

Put $\lambda^j_i=\mu^j_i t^i$ and $\lambda=(\lambda^0_i)_i$. From the
fact that ${\Bscr}_{\ge 0}$ is ample for $\mu$ (Theorem
\ref{ref-5.2.1-42}) we easily obtain that $\Cscr$ is ample for
$\lambda$.  We will now analyze the map $\Ascr\r \lambda_\ast \Cscr$.
\begin{step} The map  $\Ascr_{ii}\r
  (\lambda_\ast \Cscr)_{ii}$ is monic. If we denote its cokernel by
  $\Sscr_{ii}$ then $\Sscr_{ii}$ is locally free of
  rank one on both sides.

To see this we will show that
\begin{equation}
\pr_{1\ast}(\Ascr_{ii})\r \pr_{1\ast}(\lambda_\ast(\Cscr)_{ii})
\end{equation}
is monic and its cokernel is locally free of rank one. The corresponding statement for the right
structure is similar.

We have $\Oscr_{X_i}=\pr_{1\ast}(\Ascr_{ii})$ and
$\pr_{1\ast}(\Cscr_{ii})=\pr_{1\ast}(\lambda^0_i,\lambda^0_{i})_\ast(\Oscr_{F^0})
=\lambda^0_{i\ast}(\Oscr_{F^0})=
\lambda^0_{i\ast}\theta_\ast^{i0}(\Oscr_{F^i})=\lambda^i_{i\ast}(\Oscr_{F^i})
$.
 So we need to show that 
$\Oscr_{X_i}\r \lambda^i_{i\ast}(\Oscr_{F^i})$ is monic and that its
cokernel is locally free of rank one.

Put $B=\PP_{X_i}(\pr_{1\ast}(\Escr_i))$ and let
 $\Oscr_B(n)=\Oscr_{\PP_{X_i}(\pr_{1\ast}(\Escr_i))}(n)$. Denote the
 projection map $B\r X_i$ by $p$.
By Proposition \ref{ref-4.4.2-37} 
the map $ E^i\r
B$  is a
closed immersion.
 So the composition $F^i\r E^i\r
B$ is a closed immersion as well. We denote this composition by $v$. Now since
$\Ascr$ is of Type II it easy to see that $\dim F^i=\dim X_i$. Hence $F$
is a divisor in $B$. Generically $\Escr_0$
will be invertible over its support and hence generically $F$ will
have degree two over $X_i$. Since according to \cite[II. Ex. 7.9]{H} one has
$\Pic(B)=\Pic(X_i)\times \ZZ^x$ where $x$ is the number of connected
 components of $X$ and the factor $\ZZ^x$
corresponds to the degrees over the generic fibers it follows that 
$\Oscr_B(-F^i)=\Lscr\otimes_{\Oscr_{X_i}} \Oscr_B(-2)$ where $\Lscr\in
 \Pic(X_i)$. 

We now apply $Rp_\ast$ to the exact sequence
\begin{equation}
\label{ref-6.5-47}
0\r \Oscr_B(-F^i)\r \Oscr_B\r v_\ast\Oscr_{F^i}\r 0
\end{equation}
Using the known properties of the map $p:B\r X_i$ \cite[Ex.\
III.8.4]{H} we extract from the long exact sequence for $Rp_\ast$ a short exact 
sequence
\begin{equation}
\label{ref-6.6-48}
0\r \Oscr_{X_i}\r \lambda^i_{i\ast}(\Oscr_{F^i})\r \wedge^2
(\pr_{1\ast}\Escr_i)^\ast \otimes_{\Oscr_{X_i}}\Lscr\r 0
\end{equation}
This proves what we want.

We obtain in addition that $R^h\lambda^i_{i\ast}(\Oscr_{F^i})=0$ for
$h>0$. This may be rephrazed as the next step.
\end{step}
\begin{step}
$R^h(\lambda^0_{i},\lambda^0_{i})_\ast(\Cscr_{ii})=0$ for $h>0$. 
\end{step}
\begin{step} The map $\Ascr_{i,i+1}\r (\lambda_\ast\Cscr)_{i,i+1}$ is
  an isomorphism. 

Arguing as in Step 1 we reduce the problem to showing
that the canonical map $\pr_{1\ast}(\Escr_i)\r
\lambda^i_{i\ast}(\Oscr_{F^i}(1))$ is an isomorphism.  

Tensoring \eqref{ref-6.5-47} by $\Oscr_B(1)$ and applying
$Rp_\ast$ we obtain what we want and in addition we obtain
$R^h\lambda^i_{i\ast} (\Oscr_{F^i}(1))=0$ for $h>0$. This then
  yields the next step.
\end{step}
\begin{step}
  $R^h(\lambda^0_{i},\lambda^0_{i+1})_\ast(\Cscr_{i,i+1})=0$ for
  $h>0$. Indeed the image of $(\lambda^0_{i},\lambda^0_{i+1})$ is
  finite over $X_i$.  Thus it is sufficient to prove
  $\pr_{1\ast}R^h(\lambda^0_{i},\lambda^0_{i+1})_\ast(\Cscr_{i,i+1})=0$.
By the Leray spectral sequence this then follows from $R^h\lambda^0_{i,\ast}(\Cscr_{i,i+1})=0$ which is restatement of $R^h\lambda^i_{i\ast} (\Oscr_{F^i}(1))=0$.
\end{step}
\begin{step} Now we translate the exact sequence
  \eqref{ref-4.7-34} to our current situation. It becomes.
\[
0\r
\Cscr_{ij}\otimes_{\Oscr_{F^0}}\lambda^{0\ast}_j\pr_{1\ast}(\Qscr_j)
\r
\Cscr_{ij+1} \otimes_{\Oscr_{F^0}}
\lambda^{0\ast}_{j+1}\pr_{1\ast}(\Escr_{j+1}) \r \Cscr_{ij+2}\r 0
\]
Using Step 2 and 4 one obtains by induction that the following
sequence is exact 
\[
0\r (\lambda^0_{i},\lambda^0_j)_\ast (\Cscr_{ij})\otimes_{\Oscr_{X_j}} \Qscr_j
\r
(\lambda^0_{i},\lambda^0_{j+1})_\ast (\Cscr_{ij+1})\otimes_{\Oscr_{X_{j+1}}}
\Escr_{j+1}
\r
(\lambda^0_{i},\lambda^0_{j+2})_\ast\Cscr_{ij+2}\r 0
\]
and furthermore that
$R^h(\lambda^0_{i},\lambda^0_{j})_\ast(\Cscr_{i,j})=0$ for $h>0$.  
\end{step}
\begin{step}  The map $\Ascr_{ii+2}\r
  (\lambda_\ast \Cscr)_{ii+2}$ is an epimorphism. If we denote its kernel by
  $\Tscr_{ii+2}$  then $\Tscr_{ii+2}=\Sscr_{ii}\otimes_{\Oscr_{X_i}}\Qscr_i$.
 In particular
  $\Tscr_{ii+2}$  is locally free of
  rank one on both sides.

To prove these statements we consider the following commutative diagram with exact
rows. 
{\tiny
\begin{equation}
\label{ref-6.7-49}
\hbox{} \hskip -0.8cm\begin{CD}
0 @>>> (\lambda^0_{i},\lambda^0_{i})_\ast (\Cscr_{ii})\otimes_{\Oscr_{X_i}}
\Qscr_i
@>>>
(\lambda^0_{i},\lambda^0_{i+1})_\ast (\Cscr_{ii+1})\otimes_{\Oscr_{X_{i+1}}}
\Escr_{i+1}
@>>>
(\lambda^0_{i},\lambda^0_{i+2})_\ast\Cscr_{ii+2} @>>> 0\\
@. @AAA @A\cong AA @AAA\\
0 @>>> \Ascr_{ii}\otimes_{\Oscr_{X_i}} \Qscr_i @>>>
\Ascr_{ii+1}\otimes_{\Oscr_{X_{i+1}}} \Escr_{i+1} @>>> \Ascr_{ii+2}
@>>> 0
\end{CD}
\end{equation}
}
(the second row is the dual version of
\eqref{ref-6.1-43}). Applying the snake lemma to
\eqref{ref-6.7-49} together with Step 1 yields what we want.
\end{step} 
\begin{step}
Assume $j\ge i-1$. Then the complex
\[
0\r \Ascr_{ij}\otimes_{\Oscr_{X_j}} \Tscr_{jj+2} \r \Ascr_{ij+2}\r
(\lambda^0_{i},\lambda^0_{j+2})_\ast(\Cscr_{ij+2})\r 0
\]
is exact. 

We prove this by induction on $j$. The cases $j=i-1,i$ were covered by the
previous steps. Assume now $j\ge i+1$. We consider the following
commutative diagram with exact rows.
{\tiny
\[
\hbox{} \hskip -1.5cm\begin{CD}
@. 0 @. 0 @. 0\\
@. @AAA @AAA @AAA\\
0 @>>> (\lambda^0_{i},\lambda^0_{j})_\ast (\Cscr_{ij})\otimes_{\Oscr_{X_{j}}}
\Qscr_{j}
@>>>
(\lambda^0_{i},\lambda^0_{j+1})_\ast (\Cscr_{ij+1})\otimes_{\Oscr_{X_{j+1}}}
\Escr_{j+1}
@>>>
(\lambda^0_{i},\lambda^0_{j+2})_\ast\Cscr_{ij+2} @>>> 0\\
@. @AAA @AAA @AAA\\
0 @>>> \Ascr_{ij}\otimes_{\Oscr_{X_{j}}} \Qscr_{j} @>>>
\Ascr_{ij+1}\otimes_{\Oscr_{X_{j+1}}} \Escr_{j+1} @>>> \Ascr_{ij+2}
@>>> 0\\
@. @AAA @AAA @AAA\\
0 @>>> \Ascr_{ij-2}\otimes_{\Oscr_{X_{j-2}}} \Tscr_{j-2j}\otimes_{\Oscr_{X_{j}}}\Qscr_{j} @>>>
\Ascr_{ij-1}\otimes_{\Oscr_{X_{j-1}}}
\Tscr_{j-1,j+1}
 \otimes_{\Oscr_{X_{j+1}}}\Escr_{j+1} @>>> \Ascr_{ij}\otimes_{\Oscr_{X_{j}}}
\Tscr_{j,j+2}
@>>> 0\\
@. @AAA @AAA @AAA\\
@. 0 @. 0 @. 0 
\end{CD}
\]
}
By induction we may assume that the first two columns are exact. Hence
so is the third column.
\end{step} 
\begin{step} The canonical maps
  $\Ascr_{ij}\otimes_{\Oscr_{X_j}}{\Tscr_{jj+2}}\r \Ascr_{ij+2}$ and
  $\Tscr_{ii+2}\otimes_{X_{i+2}} \Ascr_{i+2,j+2}\r\Ascr_{ij+2}$  are
  monomorphism, and furthermore they define an isomorphism
  $\Ascr_{ij}\otimes_{\Oscr_{X_j}}{\Tscr_{jj+2}}\cong
  \Tscr_{ii+2}\otimes_{X_{i+2}} \Ascr_{i+2,j+2}$. 

To see this note that by the previous step we already know that the
first map is a monomorphism. A similar proof involving
\eqref{ref-4.8-35} shows that the second map is also a monomorphism.

Since by definition $\Tscr_{ii+2}$ goes to zero under the map $\Ascr\r
\lambda_\ast\Cscr$ we also have that  $\Tscr_{ii+2}\otimes_{X_{i+2}}
\Ascr_{i+2,j+2}$ goes to zero. Thus the image of $\Tscr_{ii+2}\otimes_{X_{i+2}}
\Ascr_{i+2,j+2}$ in $\Ascr_{i,j+2}$ lies in the image of
$\Ascr_{ij}\otimes_{\Oscr_{X_j}}{\Tscr_{jj+2}}$. By symmetry the
opposite inclusion will also hold and hence we are done.
\end{step}
\begin{step} $\Ascr$ is noetherian. 

By the previous steps we have an
  invertible ideal $\Jscr\subset \Ascr_{\ge 2}$ given by
  $\Jscr_{ij}=\Ascr_{ij-2}\otimes_{\Oscr_{X_{i-2}}}
  \Tscr_{j-2j}=\Tscr_{ii+2}\otimes_{\Oscr_{X_i}} \Ascr_{i+2j}$ in
  $\Ascr$ such that $(\Ascr/\Jscr)_{\ge 1}=\Dscr$ where $\Dscr_{\ge
  1}=\Cscr_{\ge 1}$ and $\Dscr_{ii}=\Oscr_{X_i}$. 

From the fact that $\Cscr$ is noetherian and the fact that all
$\Cscr_{ij}$ are coherent we easily obtain that $\Dscr$ is
noetherian. We may now conclude by invoking lemma \ref{ref-3.2.2-15}.
\end{step}
\section{Non-commutative deformations of Hirzebruch surfaces}
\subsection{Strongly ample sequences}
Let $\Escr$ an noetherian abelian category. For us a sequence $(O(n))_{n\in
  \ZZ}$ of objects in $\Escr$ is \emph{strongly ample} if the following conditions
hold
\begin{itemize}
\item[(A1)] For all  $\Mscr\in \Escr$ and for all $n$
there is an epimorphism $\oplus_{i=1}^t O(-n_i)\r \Mscr$ with $n_i\ge n$.
\item[(A2)] For all $\Mscr\in \Escr$ and for all $i>0$ one has
  $\Ext^i_\Escr(O(-n),\Mscr)=0$ for $n\gg 0$.
\end{itemize}
A strongly ample sequence $(O(n))_{n\in \ZZ}$ in $\Escr$ is ample in
the sense of \cite{Polishchuk1}. Hence using the methods of
\cite{AZ} or \cite{Polishchuk1} one obtains $ \Escr\cong \qgr(A) $ if
$\Escr$ is $\Hom$-finite, where $A$ is the noetherian $\ZZ$-algebra
$\oplus_{ij} \Hom_\Escr(O(-j),O(-i))$.

It would be interesting to know if a noncommutative $\PP^1$-bundle
always has an ample sequence. The next lemma is very weak but it is
sufficient for us below.
\begin{lemmas}
\label{veryweak}
Let $\Ascr$ be a non-commutative symmetric algebra over  $(X_n)_n$
(see \S\ref{ncsym}) with all $X_n$ being equal to a smooth projective
scheme $X$. Let $\Oscr_X(1)$ be an ample line bundle on $X$.  Assume
that $\Ascr_{i,i+1}$ is generated by global sections on the right for
all $i$ and that for each $m$ we have that
$\Oscr_X(-m)\otimes_{\Oscr_X} \Ascr_{mn}\otimes_{\Oscr_X} \Oscr_X(n)$
has vanishing cohomology for $n\gg 0$.

Then $\qgr(\Ascr)$ has a strongly ample sequence given by
$O(n)=\pi(\Oscr_X(n)\otimes_{\Oscr_X} e_{-n}\Ascr)$.
\end{lemmas}
\begin{proof} 
We have maps of $\gr(\Ascr)$-objects induced by the multiplication in $\Ascr$
\[
\Ascr_{i,i+1}\otimes_{\Oscr_X} e_{i+1}\Ascr\r e_i\Ascr
\]
which are surjective in degree $\ge i+1$. Since $\Ascr_{i,i+1}$ is generated by global sections
on the right these may be turned into maps
\begin{equation}
\label{aepi}
(e_{i+1}\Ascr)^{t_i}\r e_i\Ascr
\end{equation}
for certain $t_i$ which are still surjective in degree $\ge i+1$.

Let $\Mscr=\pi M$ with $M\in \gr(\Ascr)$ noetherian. Then there is
some $N$ such that $M_{\ge N}$ is generated in degree one. Hence there is some $N'$, which we will take
$\ge N$, such that there
is an epimorphism
\[
(\Oscr_X(-N') \otimes_{\Oscr_X} e_N\Ascr)^s\r M_{\ge N}
\]
which using the the maps given in \eqref{aepi} may be turned into epimorphisms
\[
(\Oscr_X(-N') \otimes_{\Oscr_X} e_{N'}\Ascr)^s\r M_{\ge N'}
\]
This implies condition (A1). 
We now compute 
\begin{align*}
\RHom_{\QGr(\Ascr)}(O(-n),\pi M)&=\RHom_X(\Oscr_X(-n),R\omega(\pi M)_n)\\&=R\Gamma(X,R\omega(\pi M)_n(n))
\end{align*}
According to \cite[Cor 3.3+proof]{Nyman5} and \cite[Lemma 3.4]{Nyman5} the map $M\r R\omega(\pi M)$ is an
isomorphism in high degree. Hence for $n\gg 0$:
\[
\RHom_{\QGr(\Ascr)}(O(-n),\pi M)=R\Gamma(X,M_n(n))
\]
Thus $\Hom_{\QGr(\Ascr)}(O(-n),-)$ has finite cohomological
dimension. To prove (A2) we may then assume that
$\Mscr=O(-m)=\pi(\Oscr_X(-m)\otimes_{\Oscr_X} e_m\Ascr)$ for $m$ large. Since in
that case
\[
M_n(n)=\Oscr_{X}(-m)\otimes_{\Oscr_X} \Ascr_{mn}\otimes_{\Oscr_X} \Oscr_X(n)
\]
we are done. 
\end{proof}

\subsection{Deformations of abelian categories}
\label{secdef}
\label{seclb}
For the convenience of the reader we will repeat the main statements
from \cite{VdBdef}. 
We first recall briefly some notions from~\cite{lowenvdb1}. Throughout $R$ will
be a commutative noetherian ring and $\mod(R)$ is its category of
finitely generated modules. \nocite{LDeD}

Let $\Cscr$ be an $R$-linear abelian category.  Then we have bifunctors $-\otimes_R-:\Cscr\times \mod(R)\r
\Cscr$, $\Hom_R(-,-):\mod(R)\times \Cscr\r \Cscr$ defined in the usual
way.  These functors may be derived in their $\mod(R)$-argument to
yield bi-delta-functors $\Tor^R_i(-,-)$, $\Ext_R^i(-,-)$.  An object $M\in \Cscr$ is
\emph{$R$-flat} if $M\otimes_R -$ is an exact functor, or equivalently if $\Tor_i^R(M,-)=0$ for
$i>0$.

By
definition (see \cite[\S3]{lowenvdb1}) $\Cscr$ is $R$-\emph{flat} if
$\Tor^R_i$ or equivalently $\Ext^i_R$ is effaceable in its
$\Cscr$-argument for $i>0$.  This implies that $\Tor^R_i$ and
$\Ext_R^i$ are universal $\partial$-functors in both arguments.

If $f:R\r S$ is a morphism of commutative noetherian rings such that
$S/R$ is finitely generated and $\Cscr$ is an $R$-linear abelian
category then $\Cscr_S$ denotes the (abelian) category of objects in
$\Cscr$ equipped with an $S$-action.  If $f$ is surjective then
$\Cscr_S$ identifies with the full subcategory of $\Cscr$ given by the
objects annihilated by $\ker f$.  The inclusion functor $\Cscr_S\r
\Cscr$ has right and left adjoints given respectively by $\Hom_R(S,-)$
and $-\otimes_R S$.

\medskip

Now assume that $J$ is an ideal in $R$ and let $\widehat{R}$ be the
$J$-adic completion of $R$. Recall that an abelian category $\Dscr$ is
said to be \emph{noetherian} if it is essentially small and all
objects are noetherian. Let $\Dscr$ be an $R$-linear noetherian
category and let $\Pro(\Dscr)$ be its category of pro-objects.  We
define $\widehat{\Dscr}$ as the full subcategory of $\Pro(\Dscr)$
consisting of objects $M$ such that $M/MJ^n\in \Dscr$ for all $n$ and
such that in addition the canonical map $M\r \invlim_n M/MJ^n$ is an
isomorphism. The category $\widehat{\Dscr}$ is $\widehat{R}$-linear.
The following is basically a reformulation of Jouanolou's results \cite{Joua}. 
\begin{propositions} (see \cite[Prop.\ \ref{ref-2.2.4-6}]{VdBdef}) $\hat{\Dscr}$ is a noetherian abelian subcategory of $\Pro(\Dscr)$. 
\end{propositions}
There is an exact functor 
\begin{equation}
\label{canonical}
\Phi:\Dscr\r \widehat{\Dscr}:M\mapsto \invlim_n M/MJ^n
\end{equation}
and we say that $\Dscr$ is complete if $\Phi$ is an equivalence of categories.
In addition we say that $\Dscr$ is \emph{formally flat} if $\Dscr_{R/J^n}$ is $R/J^n$-flat
for all $n$. 

\begin{definitions}
Assume that $\Cscr$ is an $R/J$-linear noetherian flat abelian
category. Then an $R$-deformation of $\Cscr$ is a formally
flat complete $R$-linear abelian category $\Dscr$ together with an
equivalence $\Dscr_{R/J}\cong \Cscr$. 
\end{definitions} 
In general, to simplify the notations, we will pretend
that the equivalence $\Dscr_{R/J}\cong \Cscr$ is just the identify.

Thus below we consider the case that $\Dscr$ is complete and formally flat and
$\Cscr=\Dscr_{R/J}$.   The
following definition turns out to be natural.
\begin{definitions} (see \cite[\eqref{new-1}]{VdBdef}) Assume that
  $\Escr$ is a formally flat noetherian $R$-linear abelian category. Let
  $\Escr_t$ be the full subcategory of $\Escr$ consisting of objects
  annihilated by a power of $J$. Let $M,N\in \Escr$.  Then the
  \emph{completed $\Ext$-groups} between $M$, $N$ are defined as
\begin{align*}
`\Ext_{\widehat{\Escr}}^i(M,N)&=\Ext_{\Pro(\Escr_t)}^i(M,N)
\end{align*}
\end{definitions}
An $R$-linear category $\Escr$ is said to be \emph{$\Ext$-finite} if $\Ext^i_\Escr(M,N)$ is a finitely
generated $R$-module for all $i$ and all objects $M$, $N\in \Escr$. Assuming $\Ext$-finiteness the completed $\Ext$-groups become computable. 
\begin{propositions} \cite[Prop.\ \ref{ref-2.5.3-22}]{VdBdef} Assume that $\Escr$ is a formally flat noetherian $R$-linear abelian category and that
  $\Escr_{R/J}$ is $\Ext$-finite.  Then  $`\Ext_{\widehat{\Escr}}^i(M,N)\in \mod(\hat{R})$  for $M,N\in
  \widehat{\Escr}$
and furthermore
\begin{align*}
`\Ext_{\widehat{\Escr}}^i(M,N)
&=\invlim_k \dirlim_l \Ext_{\Escr_{R/J^l}}^i(M/MJ^l,N/NJ^k)
\end{align*}
If $M$ is in addition $R$-flat then
\[
`\Ext_{\widehat{\Escr}}^i(M,N)=\invlim_k  \Ext_{\Escr_{R/J^k}}^i(M/MJ^k,N/NJ^k)
\]
\end{propositions}
The results below allow one to lift properties from $\Cscr$ to $\Dscr$.  They
follow easily from the corresponding infinitesimal results  
(\cite[Theorem A]{lowen4}, \cite[Prop.\ 6.13]{lowenvdb1}, \cite{VdBdef}).
\begin{propositions} 
\label{lifting}
  Let $M\in \Cscr$ be a flat object such that
  $\Ext^i_{\Cscr}(M,M\otimes_{R/J} J^n/J^{n+1})=0$ for $i=1,2$ and $n\ge 1$.
  Then there exists a unique  $R$-flat object (up to non-unique
  isomorphism) $\overline{M}\in \Dscr$  such that $\overline{M}/\overline{M}J\cong M$.
\end{propositions}

\begin{propositions} 
\label{vanishing}
Let $\overline{M},\overline{N}\in \Dscr$ be flat objects and put
$\overline{M}/\overline{M}J= M$, $\overline{N}/\overline{N}J= N$. Assume that for all $X$
in $\mod(R/J)$ we have $\Ext^i_\Cscr(M,N\otimes_{R/J} X)=0$ for a certain $i>0$. Then
we  have $`\Ext^i_\Dscr(\overline{M},\overline{N}\otimes_R X)=0$ for all $X\in \mod(R)$. 
\end{propositions}
\begin{propositions} 
\label{basechange}
Let $\overline{M},\overline{N}\in \Dscr$ be flat objects and put
$\overline{M}/\overline{M}J= M$, $\overline{N}/\overline{N}J=
N$. Assume that for all $X$ in $\mod(R/J)$ we have
$\Ext^1_\Cscr(M,N\otimes_{R/J} X)=0$. Then $\Hom_\Dscr(\overline{M},\overline{N})$ is
$R$-flat and furthermore 
for all $X$ in $\mod(R)$ we
have $\Hom_\Dscr(\overline{M},\overline{N}\otimes_R
X)=\Hom_\Dscr(\overline{M},\overline{N})\otimes_R X$.
\end{propositions}
Let us also mention Nakayama's lemma \cite{VdBdef}.
\begin{lemmas} \label{nakayama}  Let $M\in \Dscr$ be
    such that $MJ=0$. Then $M=0$.
\end{lemmas}
The following result is a version of ``Grothendieck's existence theorem''.
\begin{propositions} (see \cite[Prop.\ \ref{ref-4.1-37}]{VdBdef}).
\label{groht}
Assume that $R$ is complete and let $\Escr$ be an $\Ext$-finite
$R$-linear noetherian category with a strongly ample sequence
$(O(n))_n$. Then $\Escr$ is complete and furthermore if $\Escr$ is
flat then we have for $M,N\in\Escr$:
\begin{equation}
\label{extidentity}
\Ext^i_\Escr(M,N)=`\Ext^i_\Escr(M,N)
\end{equation}
\end{propositions}
The following result shows that the property of being strongly ample lifts well. 
\begin{theorems} (see \cite[Thm.\ \ref{ref-4.2-44}]{VdBdef}) \label{stronglyample} Assume
  that $R$ is complete and that $\Cscr$ is $\Ext$-finite and let
  $O(n)_n$ be a sequence of $R$-flat objects in $\Dscr$ such that
  $(O(n)/O(n)J)_n$ strongly ample. Then
\begin{enumerate}
\item $O(n)_n$ is
  strongly ample in $\Dscr$;
\item  $\Dscr$ is flat (instead of just formally flat);
\item $\Dscr$ is $\Ext$-finite as $R$-linear category. 
\end{enumerate}
\end{theorems}

\subsection{Deformations of Hirzebruch surfaces}
Below $(R,m)$ is a complete commutative local noetherian ring with
residue field $k=R/m$. Everything will now either be over $k$ or over
$R$. Although in the main part of this paper we have set up the theory over a base scheme
of finite type over a $k$ it is not difficult to see that the results remain valid over
$\Spec R$. We will use this without further comment. When
we say that something is ``compatible with base change'' we mean compatible with
the passage from $R$ to $k$. We usually abbreviate $-\otimes_R k$ by $(-)_k$. 
We also use a subscript $k$ to indicate that something is defined over $k$.

We let $X_k$ be the  Hirzebruch surface
$\PP(\Escr_k)$ with $\Escr_k=\Oscr_{\PP^1_k}\oplus \Oscr_{\PP^1_k}(h)$,
$h\ge 0$  and we let $\Dscr$ be 
an $R$-deformation of $\Cscr=\coh(X_k)$ in the sense of \S\ref{secdef}.  The rest of
this section will be devoted to proving the following result.
\begin{theorems}
\label{ref-7.4.1-70}
There exists a sheaf-bimodule $\Escr$ over
$\PP^1_R$ such that $\Dscr$ is equivalent to $\qgr(\SS(\Escr))$.
\end{theorems}

Let $t:X_k\r \PP^1_k$ be the projection map. Then we have  standard 
line bundles $\Oscr_{k}(m,n)=t^\ast
\Oscr_{\PP^1_k}(m) \otimes_{\Oscr_{X_k}} \Oscr_{X_k/\PP^1}(n)$ on $X_k$.
From the formula
\[
R\Gamma(X_k,\Oscr_{X_k}(m,n))=R\Gamma(\PP^1_k,\Oscr_{X_k}(m)
\otimes_{\Oscr_{\PP^1_k}}  Rt_\ast \Oscr_{X_k/\PP^1_k}(n))
\]
we deduce that in particular $H^i(X_k,\Oscr_{X_k}(m,n))=0$ for $i>0$ and $m,n\ge 0$.

 Since the $\Oscr_{k}(m,n)$ are exceptional in $\Cscr$ they
lift to objects $O(m,n)$ in $\Dscr$ using Proposition
\ref{lifting}.   Furthermore from the ampleness criterion in
\cite[Cor. V.2.18]{H} together with
Theorem \ref{stronglyample}(1) it follows that  $(O(n,n))_n$
is a
strongly ample sequence in $\Dscr$. By item (3) of the same theorem
we obtain
that  $\Dscr$ is $\Ext$-finite. 

We now define some $R$-linear $\ZZ$-algebras
\[
C_{n}=\bigoplus_{j\ge i}\Hom(O(-j,-n),O(-i,-n))
\]
as well as $C_{m}-C_{n}$-bimodules for $n\ge m$:
\[
A_{mn}=\bigoplus_{j\ge i}\Hom(O(-j,-n),O(-i,-m))
\]
From Proposition \ref{basechange} it follows that $C_n$ and $A_{mn}$ are $R$-flat
and compatible with base change. Hence
\begin{align}
\label{Cdesc} C_{n,k}&=\bigoplus_{j\ge i} \Gamma(\PP^1_k,\Oscr_{\PP^1_k}(j-i))\\
\label{Adesc} A_{mn,k}&=\bigoplus_{j\ge i}\Hom(S^{n-m}\Escr_k(j-i))
\end{align}

We can now look for some properties of $C_{n,k}$ that lift to $C_{n}$ 
(see \cite[\S8.3]{VdB26} for a more elaborate example of how this is done).
\begin{itemize}
\item[(P1)] 
\[
\rk C_{n,ij}=
\begin{cases}
j-i+1&\text{if  $j\ge i$}\\
0&\text{otherwise}
\end{cases}
\]
\item[(P2)]
Define $V_{n,i}=C_{n,i,i+1}$. Then $C_{n}$ is generated by the $(V_{n,i})_i$.
\item[(P3)]
Put $K_{n,i}=\ker(V_{n,i}\otimes V_{n,i+1}\r C_{n,i,i+2})$.  Then the relations
between the $V_{n,i}$ in $C_{n}$ are generated by the $K_{n,i}$.
\item[(P4)]
Rank counting
reveals that  $\rk K_{n,i}=1$.  The $R$-module $K_{n,i}$ is generated by a non-degenerate
tensor $r_{n,i}$ in $V_{n,i}\otimes_R V_{n,i+1}$.
\end{itemize}
Using these properties it is now easy to describe $C_{n}$. After
choosing suitable bases $x_i,y_i$ in $V_{n,i}$ we may assume that
$r_{i}=y_i x_{i+1}-x_i y_{i+1}$. Thus all $C_{n}$ are in fact
isomorphic to $\check{S}$ (see \eqref{ref-4.4-26}) where $S$ is the graded
algebra $R[x,y]$. In particular $\qgr(C_n)\cong \coh(\PP^1_R)$ for all $n$.

It also follows that after suitable reindexing $A_{mn}$ becomes in a
natural way a bigraded $S\otimes_R S$-module which we denote by $A'_{mn}$.  We think
of $A'_{mn}$ as an $S$-$S$-bimodule with independent left and right grading. 
The
required reindexing is given by
\[
A'_{mn;ij}=A_{mn;-i,j}
\]
Here $x,y$ act as $x_{i-1},y_{i-1}$ on the left and as $x_j,y_j$ on the right. 

The following diagram is commutative
\[
\xymatrix{
  \gr(C_m)\ar@{=}[d]\ar[rr]^{-\otimes_{C_m}A_{mn}}&&\gr(C_n)\ar@{=}[d]\\
\gr(S)\ar[rr]^{(-\otimes_S A'_{mn})_{0,-}}&&\gr(S)
}
\]
Here by $(-)_{0,-}$ we mean taking the part of degree zero for the
left grading.  

\medskip

Let $\Ascr_{mn}$ be the quasi-coherent
$\Oscr_{\PP^1_R}\boxtimes\Oscr_{\PP^1_R}$-module associated to
$A'_{mn;ij}$.
\begin{lemmas} 
\label{explicitA}
$\Ascr_{mn}$ is locally free on the left and right of rank $n-m+1$. In addition
\[
\Ascr_{mn,k}\cong \delta_\ast S^{n-m}\Escr_k
\]
where $\delta:\PP^1_k\r \PP^1_k\times \PP^1_k$ is the diagonal
embedding.  
\end{lemmas}
\begin{proof}
We first observe that $\Ascr_{mn}$ is in fact coherent. 
To this end it is sufficient
to show that the diagonal submodule $\bigoplus_{i} A'_{mn;ii}$ is a finitely
generated $\bigoplus_{i} S_i\otimes_R S_i$-module. This may be verified after tensoring
with $k$.

From \eqref{Adesc} one obtains
\begin{equation}
\label{wasalign}
\begin{split}
A_{mn}'\otimes_R k&=\bigoplus_{j+i\ge 0}\Gamma(X_k,O_k(j+i,n-m))\\
&=\bigoplus_{j+i\ge 0}\Gamma(\PP^1_k,S^{n-m}\Escr_k(j+i))
\end{split}
\end{equation}
Thus
\begin{equation}
\label{rhs37}
\bigoplus_i A_{mn;ii}'\otimes_R k=\bigoplus_{i\ge 0}\Gamma(\PP^1_k,S^{n-m}\Escr_k(2i))
\end{equation}
The righthand side of \eqref{rhs37} is the graded-$\bigoplus_i
S_{i,k}\otimes_k S_{i,k}$-module associated to the coherent
$\Oscr_{\PP^1_k\times \PP^1_k}$-module $\delta_\ast S^{n-m}\Escr_k$ (for the ample line
bundle given by $\Oscr_k(1,1)$). Hence this graded module is finitely generated.

From the computation in the previous
paragraph we also learn that $\Ascr_{m,n}\otimes_R k$ is indeed given by the sheaf
$S^{n-m}\Escr_k$ supported on the diagonal. 

We claim that the  support of $\Ascr_{mn}$  is finite
over both factors of $\PP^1_R\times\PP^1_R$. Again it is clearly
sufficient to check this over $k$ but then it follows from the
explicit form of $\Ascr_{m,n}\otimes_R k$ given above.

As indicated above $A_{mn}$ is flat over $R$. Hence the same is true for $\Ascr_{mn}$.
Since $\Ascr_{mn}\otimes_R k$ is locally free over both factors it follows from
Lemma \ref{ref-3.1.5-4} that $\Ascr_{m,n}$ is
locally free on the left and on the right. By tensoring with $k$ 
we deduce that the left and right rank of $\Ascr_{m,n}$ are equal to
$n-m+1$.
\end{proof}
\begin{lemmas} The functor $-\otimes_{C_m}A_{mn}$ sends $\gr(C_m)$ to $\gr(C_n)$.
\end{lemmas}
\begin{proof} It is sufficient to prove that for every $i$ we have that $e_i A_{mn}$ 
lies in $\gr(C_n)$. Since $e_i A_{mn}$ is a finitely generated $R$-module in
every degree we may prove this after specialization.

We compute
\begin{align*}
e_i A_{mn,k}=\bigoplus_{j\ge i} \Gamma(\PP^1_k,S^{n-m}\Escr_k(j-i))
\end{align*}
Thus $e_i A_{mn,k}$ is up to finite length modules the graded
$S$-module associated to the coherent $\PP^1$-module $S^{n-m}\Escr_k(-i)$. Hence
it is finitely generated. 
\end{proof}
\begin{lemmas}
\label{diaglemma}
There is a commutative diagram
\begin{equation}
\label{comparison}
\xymatrix{
  \gr(C_m)\ar[rr]^-{-\otimes_{C_m}A_{mn}}\ar[d]_\pi&&\gr(C_n)\ar[d]^\pi\\
\coh(\PP^1_R)\ar[rr]_{-\otimes_{\PP^1_R}{\Ascr_{mn}}}&&\coh(\PP^1_R)
}
\end{equation}
\end{lemmas}
\begin{proof}
We first have to construct a natural transformation
\[
\xymatrix{
  \gr(C_m)\ar[rr]^-{-\otimes_{C_m}A_{mn}}\ar[dd]_\pi&&\gr(C_n)\ar[dd]^\pi\\
&&\ar@{=>}[dl]\\
\coh(\PP^1_R)\ar[rr]_{-\otimes_{\PP^1_R}{\Ascr_{mn}}}&&\coh(\PP^1_R)
}
\]
Taking into account the equivalences $\gr(C_m)=\gr(S)$ this diagram may be rewritten
as 
\begin{equation}
\label{rewrittendiag}
\xymatrix{
  \gr(S)\ar[rrrr]^-{(-\otimes_{S}A'_{mn})_{0,-}}\ar[dd]_\pi&&&&\gr(S)\ar[dd]^\pi\\
&&&&\ar@{=>}[dl]\\
\coh(\PP^1_R)\ar[rrrr]_{\pi M\mapsto \pi([\omega_1\pi_1(M\otimes_S A'_{mn})]_{0,-})}&&&&\coh(\PP^1_R)
}
\end{equation}
Here $\omega_1$ is $\omega$ applied to the left
grading and similarly for $\pi$. The natural transformation is now
obtained by functoriality from the canonical map
\[
M\otimes_{S}A'_{mn}\r \pi_1\omega_{1}(M\otimes_S A'_{mn})
\]
We claim this natural transformation is an isomorphism.
Both branches of the diagram \eqref{comparison} represent right exact functors so it is sufficient
to consider the value on the projective generators $S_k(i)$ of $\gr(S_k)$. This verification
may be done after specialization.

We find
\begin{align*}
S_k(i)\otimes_{S_k} A'_{mn,k}&=\bigoplus_{pq} A'_{mn;p+i,q,k}\\
&=\bigoplus_{p+q+i\ge 0} \Gamma(\PP^1_k,S^{n-m}\Escr_k(q+p+i))
\end{align*}
where have used \eqref{wasalign}.
An easy verification shows that
\[
\pi_1\omega_{1}(S_k(i)\otimes_{S_k} A'_{mn,k})=\bigoplus_{p,q} \Gamma(\PP^1_k,S^{n-m}\Escr_k(q+p+i))
\]
from which we deduce that
\[
(\pi_1\omega_{1}(S_k(i)\otimes_{S_k} A'_{mn,k})/(S_k(i)\otimes_{S_k} A'_{mn,k}))_{l,-}
\]
is finite dimensional for any $l$. This implies that the natural transformation
in \eqref{rewrittendiag} is in fact a natural isomorphism. 
\end{proof}
The natural morphism $A_{mn}\otimes_{C_n} A_{nt}\r A_{mt}$ induces via diagram \eqref{comparison}
a natural transformation of functors
\[
-\otimes_{\PP^1_R}(\Ascr_{mn}\otimes_{\PP^1_R}\Ascr_{nt})\r -\otimes_{\PP^1_R}\Ascr_{mt}
\]
Using Lemma \ref{ref-3.1.1-3} one obtains from this a morphism of bimodules
\begin{equation}
\label{aprod}
\Ascr_{mn}\otimes_{\PP^1_R}\Ascr_{nt}\r \Ascr_{mt}
\end{equation}
Using a similar argument one shows that this morphism of bimodules satisfies the
associativity axiom and hence produces a sheaf-$\ZZ$-algebra on $\PP^1_R$ given
by
\[
\Ascr=\bigoplus_{n\ge m}\Ascr_{mn}
\]
From the fact that $A_{mm}=C_m$ on easily obtains $\Ascr_{mm}=\Oscr_{\PP^1_R}$. 
Explicitating the proof of Lemma \ref{ref-3.1.1-3} one obtains that over $k$ \eqref{aprod} is
given by the canonical maps
\[
S^{n-m}\Escr_k\otimes_{\PP^1_k} S^{t-n} \Escr_k\r S^{t-m}\Escr_k
\]
Therefore by a suitable version of Nakayama's lemma we deduce that
\eqref{aprod} is an epimorphism and hence $\Ascr$ is generated by $\Escr_n\overset{\text{def}}{=}\Ascr_{n,n+1}$.

Let $Q_n$ be the kernel of
$\Ascr_{n,n+1}\otimes_{{\PP^1_R}}\Ascr_{n+1,n+2}\r
\Ascr_{n,n+2}$. We claim that this  kernel is non-degenerate in
$\Ascr_{n,n+1}\otimes_{{\PP^1_R}}\Ascr_{n+1,n+2}$. 

In Lemma \ref{ref-3.1.10-12} we have shown that dualizing of bimodules is
compatible with base change. 
From this it easily follows that  it is
sufficient to check the non-degenerateness of $\Qscr_n$ over $k$ where
it is obvious.  

Now let $\Ascr'$ be the $\ZZ$-algebra generated by the
$\Escr_{n,n+1}$ subject to the relations given by the $\Qscr_n$.
By construction there is  a surjective map  $\Ascr'\r \Ascr$. Since $\Ascr'$ and $\Ascr$ are
locally free in each degree and have the same rank it follows that
this surjective map must actually be an isomorphism.

So summarizing we have shown the following:
\begin{lemmas} $\Ascr$ is a non-commutative symmetric algebra over $\PP^1_R$.
\end{lemmas}
By \S\ref{ref-4.1-21} it follows that $\gr(\Ascr)\cong 
\gr(\SS(\Escr))$ with $\Escr=\Ascr_{01}$ and this equivalence preserves right bounded modules. Hence to finish the proof of Theorem \ref{ref-7.4.1-70} it
  is sufficient to show  
 that $\qgr(\Ascr)\cong \Dscr$.

Put 
\[
C=\bigoplus_{\begin{smallmatrix}(i,m),(j,n)\\ i\le j\\m\le n\end{smallmatrix}}
\Hom(O(-j,-n),O(-i,-m))
\]
Then $C$ is a $\ZZ^2$-algebra and we have an exact functor 
\[
\Sigma:\Gr(C)\r \Gr(\Ascr)
\]
which is defined as follows. Let $M\in \Gr(C)$. Then $M_n\overset{\text{def}}{=}M_{-,n}$
is a right $C_n$-module. Furthermore the right action of $C$ on $M$ induces
maps
\begin{equation}
\label{nat37}
M_m\otimes_{C_m}A_{mn}\r M_n
\end{equation}
Put $\Mscr_n=\pi(M_n)\in \Qch(\PP^1_R)$. Thanks to Lemma \ref{diaglemma} the
maps \eqref{nat37} become maps
\[
\Mscr_m\otimes_{\PP^1_R}\Ascr_{mn}\r \Mscr_n
\]
and one checks that $\Sigma\Mscr\overset{\text{def}}{=}\bigoplus_n
\Mscr_n$ defines an object in $\Gr(\Ascr)$. Put $\sigma M=\pi\Sigma
M\in\QGr(\Ascr)$ (where here $\pi$ is the quotient functor $\Gr(\Ascr)\r \QGr(\Ascr)$).

We claim that $\Sigma$ sends finitely generated objects in $\Gr(C)$ to objects in
$\gr(\Ascr)$. It suffices to prove this for the projective generators $e_{im}C$.

We have for $n\ge m$
\[
e_{im}C_{-,n}=\bigoplus_{j\ge i} \Hom_\Dscr(O(-j,-n),O(-i,-m))
\]
Hence we have to prove that the righthand side is a finitely generated $C_n$-module.
Since the summands $\Hom_\Dscr(O(-j,-n),O(-i,-m))$ are all finitely generated
$R$-modules we may do this after specialization. We get
\[
e_{im}C_{-,n}\otimes_R k=\bigoplus_{j\ge i}\Gamma(\PP^1_k,S^{n-m}\Escr_k(j-i))
\]
which is indeed finitely generated. For reference below we note that
from this computation we also get
\[
\pi(e_{im}C_{-,n}\otimes_R k)=S^{n-m}\Escr_k(-i)
\]
(where here $\pi$ is the quotient functor $\Gr(C_{n,k})\r
\QGr(C_{n,k})\cong \Qch(\PP^1_k)$) and thus
\[
\Sigma(e_{im}C\otimes_R k)=\bigoplus_{n\ge m} S^{n-m}\Escr_k(-i)
\]
so that finally we get
\begin{equation}
\label{onemorefinally}
\sigma(e_{im} C\otimes_R k)=O_k(-i,-m)
\end{equation}

Since
\[
\Hom_C(e_{jn}C,e_{im}C)=e_{im}Ce_{jn}=C_{(i,m)(j,n)}=\Hom_\Dscr(O(-j,-n),O(-i,-m))
\]
functoriality yields a morphism of $R$-modules (for $j\ge i$, $n\ge m$)
\[
\Hom_\Dscr(O(-j,-n),O(-i,-m))\r \Hom_{\qgr(\Ascr)}(\sigma(e_{jn}C),\sigma(e_{im}C))
\]
The left hand side is $R$-flat and commutes with base change as indicated above. 
We claim that this true for the right hand side as well.
\begin{lemmas} 
\label{deformationA} $\qgr(\Ascr)$ is a deformation of $\qgr(\Ascr)_k=\qgr(\Ascr_k)=\qgr(S\Escr_k)=\coh(\PP^1_k)$.
\end{lemmas}
\begin{proof} According to \cite{Nyman5} $\qgr(\Ascr)$ is $\Ext$-finite.
Therefore, according to Proposition \ref{groht} it is sufficient to prove that
  $\qgr(\Ascr)$ has a strongly ample sequence. To this end we verify the
  conditions for Lemma \ref{veryweak}.  It is standard that these
  conditions lift from $k$ to $R$ and hence we may check them over
  $k$. Over~$k$ they follow from the explicit description of
  $\Ascr_{mn,k}$ given in Lemma \ref{explicitA}.
\end{proof}
\begin{lemmas}
The $R$-module
\[
 \Hom_{\qgr(\Ascr)}(\sigma(e_{jn}C),\sigma(e_{im}C))
\]
for $i\le j$ and $m\le n$ is flat and compatible with base
change. Furthermore the canonical map
\begin{equation}
\label{Aisomorphism}
\Hom_\Dscr(O(-j,-n),O(-i,-m))\r \Hom_{\qgr(\Ascr)}(\sigma(e_{jn}C),\sigma(e_{im}C))
\end{equation}
constructed above is an isomorphism. 
\end{lemmas}
\begin{proof} We first discuss the first statement. Given
  Lemma \ref{deformationA} it is sufficient to check that
  $\sigma(e_{jn}C)\otimes_R k$ satisfies the conditions of Proposition
  \ref{basechange}. It is easy to see that $\sigma(e_{jn}C)$ is compatible
  with base change and is $R$-flat. One may then invoke the explicit
  description of $\sigma(e_{jn}C\otimes_R k)$ given in
  \eqref{onemorefinally}.

  To prove the last statement we note that is is true over $k$ by
  \eqref{onemorefinally}. We may then invoke Nakayama's lemma for $R$ (given that
everything is compatible with base change as we have shown above). 
\end{proof}
\begin{proof}[Proof of Theorem \ref{ref-7.4.1-70}]
  Given our preparatory work it is sufficient to prove $\Dscr\cong
  \qgr(\Ascr)$.  By Theorem \ref{stronglyample} we obtain that
  $(O(n,n))_n$ is an ample sequence in~$\Dscr$.  Given
  \eqref{Aisomorphism} and the $\ZZ$-algebra version of the
  Artin-Zhang theorem \cite{AZ} it is sufficient to prove that
  $(\sigma(e_{-n,-n}C))_n$ forms a strongly ample sequence in
  $\qgr(\Ascr)$. Using Lemma \ref{deformationA} together with Theorem
  \ref{stronglyample} this may be checked over $k$. Then we invoke
  again the explicit description of $\sigma(e_{-n,-n}C\otimes_R k)$
  given in \eqref{onemorefinally}.
\end{proof}


\def\cprime{$'$} \def\cprime{$'$} \def\cprime{$'$}
\ifx\undefined\bysame
\newcommand{\bysame}{\leavevmode\hbox to3em{\hrulefill}\,}
\fi

\end{document}